\documentclass[10pt]{article}
\usepackage[pass,letterpaper, total={4.5in, 7.125in},bindingoffset=0in,left=0in,right=0in,top=0in,bottom=0in,footskip=0in]{geometry}

\usepackage{anyfontsize}
\usepackage{lmodern}
\usepackage{enumitem}
\usepackage{amsthm,amsmath,amssymb}
\usepackage{newtxtext}
\usepackage{newtxmath}

\usepackage{graphicx}
\usepackage{tikz}
\usepackage{tkz-berge}
\usetikzlibrary{shapes.geometric}
\newcommand{\tikzLongrightarrow}{%
 \mathbin{\tikz{\draw[arrows={-latex},line width=1.2pt,double=white] (0,0) -- (3em,0);}}
}
\usepackage[colorlinks=true,citecolor=black,linkcolor=black,urlcolor=blue]{hyperref}

\renewcommand{\mod}[1]{{\rm ~\mkern-1mu (mod~}#1)}
\usepackage{latexsym}
\usepackage{amssymb}

\newlength\triplesep
 \newlength\triplelinewidth 
 \tikzset{triple/.style={line width=\triplelinewidth,green, preaction={ preaction={draw,line width=2\triplesep+3\triplelinewidth,green}, draw,line width=2\triplesep+\triplelinewidth,white} } }
 \tikzset{quple/.style={line width=.5pt, preaction={draw,double,double distance=2.5pt,line width=.5pt,gray}, draw,double, line width=.5pt,double distance=.5pt,gray} } 
\tikzstyle{vertex}=[circle,fill=green!25, draw=black,minimum size=26pt,inner sep=0pt]
\tikzstyle{nvertex}=[circle,fill=green!25, draw=black,minimum size=0pt,inner sep=0pt]
\tikzstyle{edge} = [draw,thick,-]
\tikzstyle{ledge} = [draw,-]
\tikzstyle{dedge} = [draw,thick, purple, double,double distance=1pt,-]
\tikzstyle{tredge} = [draw,thick, triple,-]
\tikzstyle{math}=[]
\tikzstyle{tedge} = [draw,line width=0.6mm,color=blue,-]
\tikzstyle{enode}=[circle, draw, fill=black!50, inner sep=0pt, minimum width=13pt]
\tikzstyle{snode}=[circle, draw, fill=black!50, inner sep=0pt, minimum width=6pt]\tikzstyle{pnode}=[regular polygon,regular polygon sides=5, draw, fill=black!80, inner sep=0pt, minimum width=17pt]
\usetikzlibrary{positioning}

\begin{document}

\newtheorem{mydef}{Definition}
\newtheorem{lem}{Lemma}
\newtheorem{cor}{Corollary}
\newtheorem{thm}{Theorem}
\newtheorem{obs}{Observation}
\newtheorem{lemma}{Lemma}
\newtheorem{sublemma}{Lemma}[lemma]
\renewcommand{\emptyset}{\varnothing}
\newcommand{\dtri}{{\rm diamond}}
\newcommand{\W}{\mathbb W}
\newcommand{\nj}[1]{\langle#1\rangle}
\newcommand{\contrv}{{/}}

%\date{21$^{\rm st}$ July, 2017\\{\small Mathematics Subject Classifications: 05C31, 05C69}}

\title{{\textbf {Quintic graphs with every edge in a triangle}}}
\author{
James Preen\\
\small Mathematics,\\[-0.8ex]
\small Cape Breton University,\\[-0.8ex]
\small Sydney, Nova Scotia, B1P 6L2, Canada.\\
\small\tt james\_preen@capebretonu.ca\\
}
\date{}
\maketitle

\begin{abstract}
We characterise the quintic (i.e. 5-regular) multigraphs with the property
that every edge lies in a triangle. Such a graph is either from a set of small graphs or is formed by adding a perfect matching to a line graph of a cubic graph as double edges,
or can be reduced by a sequence of %local 
operations to one of these graphs. 
\end{abstract}

\section{{Introduction}}\label{s:intro}
 In this paper a triangle in a graph will be defined as a set of three distinct vertices with an edge between each pair of vertices. The open neighbourhood of a set $S$ of vertices is the union of the sets of vertices adjacent to a each vertex in $S$, with all vertices from $S$ removed. 
% A configuration will be an induced subgraph of a graph.
For other basic graph theory definitions, please see \cite{bo:rd}.

We are interested in graphs such that every edge is in at least one triangle, we refer to this as {\em the triangle property}. 
By our definition of a triangle multiple edges are permitted in graphs with the triangle property but loops cannot be in a triangle since they must repeat a vertex. 
In \cite{ar:PR} a similar characterisation was given for 4-regular graphs, and the 5-regular case was mentioned as the question that originally motivated their paper.

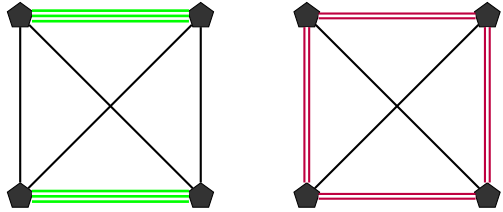
\begin{figure}[ht]
\begin{center}
\begin{tikzpicture}[scale=0.6, transform shape]
   \foreach \pos/\name in {{(0,0)/a}, {(4,0)/b}, {(4,4)/c},{(0,4)/d}}
        \node[pnode] (\name) at \pos {};
   \foreach \source/\dest in {b/c,b/d,d/a,c/a}
        \path[edge] (\source) -- (\dest);
    \foreach \source/\dest in {a/b,c/d}
        \path[tredge] (\source) -- (\dest);
%\node[math] () at (-1,2) {$H$:};
\end{tikzpicture}
~~~~~~~~~~
\begin{tikzpicture}[scale=0.6, transform shape]
   \foreach \pos/\name in {{(0,0)/a}, {(4,0)/b}, {(4,4)/c},{(0,4)/d}}
        \node[pnode] (\name) at \pos {};
   \foreach \source/\dest in {b/c,c/d,d/a,b/a}
        \path[dedge] (\source) -- (\dest);
    \foreach \source/\dest in {a/c,b/d}
        \path[edge] (\source) -- (\dest);
%\node[math] () at (-1,2) {$H$:};
\end{tikzpicture}
\end{center}
\caption{The smallest quintic graphs with the triangle property}\label{f:fourv}
\end{figure}

Any 5-regular graph must have an even number of vertices and a triangle requires three distinct vertices, so the smallest quintic graphs with the triangle property have four vertices, as shown in figure \ref{f:fourv}. If there were a pair of vertices which were not joined, then they would each require three extra edges from them to other vertices, but the degrees of their neighbours are at least 3 and so only 4 edges can be added from the unjoined vertices as we are in a quintic graph.
Using the complete graph $K_4$ as a base, we need to add four edges with a maximum of two added to any vertex, and that is either two double edges or a 4 cycle.

We will use the graphical conventions that a vertex with all five of its edges going to vertices shown in the figure is a pentagon. Multiple edges will be shown by parallel edges and colour is used to more easily distinguish them.
$K_6$ is the smallest simple quintic graph, and it has the triangle property, in fact every edge of it is in four triangles.

The same fundamental lemma as in the 4-regular paper will serve as our basis, and a corollary from \cite{ar:jp} that follows from it will give us more information about the structure of quintic graphs with the triangle property.

\begin{lem}\label{l:triv} \cite{ar:PR}
A graph $G$ has the triangle property if and only if, for every vertex $v \in V (G)$,
the graph induced by the neighbourhood of $v$ contains no isolated vertices.
\end{lem}

\begin{cor}\label{c:triv} \cite{ar:jp}
In a 5-regular simple graph with $2m$ vertices which has the triangle property, there are at least $m$ edges which are in two triangles or more.
\end{cor}

This implies that at every vertex of a simple graph there are two triangles which share an edge. The induced subgraph underlying this, the {\em diamond}, which is $K_4$ with an edge removed, will be important in this paper. In the case of $G$ having multiple edges it isn't necessary to have a diamond at every vertex, but if there isn't, there must be a multiple edge at that vertex. 

%\begin{proof}
%If $v$ is a vertex in a 5-regular graph $G$ with the triangle property then we can use Lemma \ref{l:triv}; if $v$ is adjacent to five distinct vertices in $G$ then the induced neighbourhood graph of $v$ has at least 3 edges, and hence a vertex $w$ of degree at least 2, which gives rise to two triangles containing $v$ and $w$. Similarly, if $e$ is a multiple edge incident with $v$ then there must be at least one other edge incident with $v$ in a triangle with $e$ which is therefore in at least two triangles.
%Since at all $2m$ vertices there is at least one edge in two triangles or more, there are at least $m$ edges in $G$ which are in at least two triangles.
%\end{proof}

A common deduction used throughout involves vertices with four known neighbours. 
If certain subgraphs are forbidden, or if the neighbour already has 5 known neighbours, there are fewer possibilities for it to be joined to.

\begin{obs}\label{o:four}
If $v$ is a vertex of degree 4 in a configuration from a 5-regular graph with the triangle property it is necessary that the fifth edge from $v$ must either be to a vertex which is adjacent to a neighbour or a multiple edge to an existing neighbour. 
\end{obs}

During this paper we will proceed by gradually showing various reductions which will make a smaller quintic graph with the triangle property, thus allowing us to then restrict further which induced subgraphs can exist in the graph. This will continue until we only have a small number of induced subgraphs which are irreducible using our methods. We can use some of these to reduce others, and then show there are only a few ways to combine these subgraphs to make an irreducible graph, giving us the base from which all quintic graphs with the triangle property can be constructed by reversals of the reductions.

\subsection{{Cut vertices in quintic graphs with the triangle property}}\label{s:lowconn}

\begin{figure}[ht]
\begin{center}
%\begin{tabular}{llll}
\begin{tikzpicture}[scale=0.4, transform shape]
\begin{scope}[shift={(-4,0)}]
\draw (1.3,0) ellipse (0.65cm and 1.3cm);
   \foreach \pos/\name in {{(0,0)/a}, {(1,.7)/b}, {(1,-.7)/d}}
        \node[enode] (\name) at \pos {};
   \foreach \pos/\name in {{(-.7,-.2)/a1},{(-.7,+.2)/a2}}
        \node[nvertex] (\name) at \pos {};
    % First we draw the vertices
    % Connect vertices with edges and draw weights
    \foreach \source/\dest in {a/d}
        \path[edge] (\source) -- (\dest);
            \foreach \source/\dest in {a/b}
        \path[dedge] (\source) -- (\dest);
\foreach \source/\dest in {a/a1,a/a2}
        \path[edge] (\source) -- (\dest);
\node[math] () at (1.4,0) {\LARGE $H$};
\node[math] () at (2.5,0) {\LARGE or};
\end{scope}
\draw (1.3,0) ellipse (0.65cm and 1.3cm);
   \foreach \pos/\name in {{(0,0)/a}, {(1,1)/b}, {(1,0)/c},{(1,-1)/d}}
        \node[enode] (\name) at \pos {};
   \foreach \pos/\name in {{(-.7,-.2)/a1},{(-.7,+.2)/a2}}
        \node[nvertex] (\name) at \pos {};
    % First we draw the vertices
    % Connect vertices with edges and draw weights
    \foreach \source/\dest in {a/b,a/c,a/d}
        \path[edge] (\source) -- (\dest);
\foreach \source/\dest in {a/a1,a/a2}
        \path[edge] (\source) -- (\dest);
\node[math] () at (1.6,0) {\LARGE $H$};
\node[math] () at (2.8,0) {$\tikzLongrightarrow$};
%\end{tikzpicture}
%&
%\begin{tikzpicture}[scale=0.8, transform shape]
%\begin{scope}[shift={(0,5)}]
   \foreach \pos/\name in {{(4.5,0)/a}}
       \node[enode] (\name) at \pos {};
   \foreach \pos/\name in {{(6.3,1)/b}, {(5.5,0)/c},{(6.3,-1)/d}}
       \node[pnode] (\name) at \pos {};
   \foreach \pos/\name in {{(3.8,-.2)/a1},{(3.8,.2)/a2}}
        \node[nvertex] (\name) at \pos {};
    % First we draw the vertices
    % Connect vertices with edges and draw weights
    \foreach \source/\dest in {a/b,a/c,a/d}
        \path[edge] (\source) -- (\dest);
\foreach \source/\dest in {a/a1,a/a2}
        \path[edge] (\source) -- (\dest);
\foreach \source/\dest in {b/c,c/d,b/d}
        \path[dedge] (\source) -- (\dest);
\node[math] () at (7.2,0) {\LARGE or};
   \foreach \pos/\name in {{(8.5,0)/a}}
       \node[enode] (\name) at \pos {};
   \foreach \pos/\name in {{(9.5,1)/b}, {(10.5,0)/c},{(9.5,-1)/d}}
       \node[pnode] (\name) at \pos {};
   \foreach \pos/\name in {{(7.8,-.2)/a1},{(7.8,.2)/a2}}
        \node[nvertex] (\name) at \pos {};
    % First we draw the vertices
    % Connect vertices with edges and draw weights
    \foreach \source/\dest in {a/d,b/d}
        \path[edge] (\source) -- (\dest);
\foreach \source/\dest in {a/a1,a/a2}
        \path[edge] (\source) -- (\dest);
\foreach \source/\dest in {a/b,c/b}
        \path[dedge] (\source) -- (\dest);
\foreach \source/\dest in {c/d}
        \path[tredge] (\source) -- (\dest);

%\end{scope}
\end{tikzpicture}
%&
~,~~~
\begin{tikzpicture}[scale=0.4, transform shape,xshift=3cm]
\draw (1.3,0) ellipse (0.65cm and 1.3cm);
   \foreach \pos/\name in {{(0,0)/a}, {(1,.8)/b}, {(1,-.8)/c}}
        \node[enode] (\name) at \pos {};
   \foreach \pos/\name in {{(-.7,-.2)/a1},{(-.8,0)/a2},{(-.7,+.2)/a3}}
        \node[nvertex] (\name) at \pos {};
    % First we draw the vertices
    % Connect vertices with edges and draw weights
    \foreach \source/\dest in {a/b,a/c}
        \path[edge] (\source) -- (\dest);
\foreach \source/\dest in {a/a1,a/a2,a/a3}
        \path[edge] (\source) -- (\dest);
\node[math] () at (1.4,0) {\LARGE $H$};
\node[math] () at (2.95,0) {$\tikzLongrightarrow$};
 %\end{tikzpicture}
%&
%\begin{tikzpicture}[scale=0.8, transform shape]
   \foreach \pos/\name in {{(5,0)/a}}
        \node[enode] (\name) at \pos {};
   \foreach \pos/\name in {{(6,.8)/b}, {(6,-.8)/c}}
        \node[pnode] (\name) at \pos {};
   \foreach \pos/\name in {{(4.3,-.2)/a1},{(4.2,0)/a2},{(4.3,+.2)/a3}}
        \node[nvertex] (\name) at \pos {};
    % First we draw the vertices
    % Connect vertices with edges and draw weights
    \foreach \source/\dest in {a/b,a/c}
        \path[edge] (\source) -- (\dest);
\foreach \source/\dest in {a/a1,a/a2,a/a3}
        \path[edge] (\source) -- (\dest);
\foreach \source/\dest in {b/c}
        \path[quple] (\source) -- (\dest);
\end{tikzpicture}
%\end{tabular}
\end{center}
\caption{Cut vertex reductions}\label{f:cutv}
\end{figure}
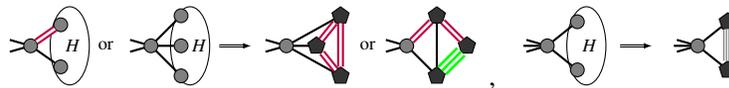

Since the triangle property is necessary for all components of a disconnected graph with the property, we can assume that all of our graphs are connected. We may choose to disconnect them in the process of reducing them, though.
As a first example of reductions, let us consider both kinds of cut-vertices possible in a 5-regular graph with the triangle property, as shown in figure \ref{f:cutv}, where 
an unknown subgraph is shown as an ellipse, half edges indicate unknown neighbours outside those shown, and the arrow gives the reduced configuration(s).
We cannot have a cut-edge since that edge could not be in a triangle. 

Thus, if there are more vertices than edges from a cut-vertex to a subgraph then we can reduce to a graph with fewer vertices. If there is a cut-vertex with a multiple edge adjacent then it must be a double edge and part of an edge cut of cardinality 3 and can again be reduced using either of the left hand reductions since $H$ will have at least four vertices. 

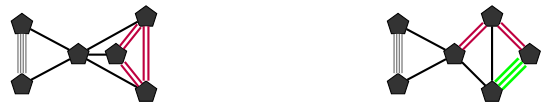
\begin{figure}[ht]
\begin{center}
\begin{tikzpicture}[scale=0.5, transform shape]
    \foreach \pos/\name in {{(4.5,0)/a},{(6.3,1)/b}, {(5.5,0)/c},{(6.3,-1)/d},{(3,.8)/y}, {(3,-.8)/z}}
       \node[pnode] (\name) at \pos {};
  \foreach \source/\dest in {a/b,a/c,a/d}
        \path[edge] (\source) -- (\dest);
\foreach \source/\dest in {b/c,c/d,b/d}
        \path[dedge] (\source) -- (\dest);
   \foreach \source/\dest in {a/y,a/z}
        \path[edge] (\source) -- (\dest);
\foreach \source/\dest in {y/z}
        \path[quple] (\source) -- (\dest);
        
\begin{scope}[shift={(10,0)}]
   \foreach \pos/\name in {{(4.5,0)/a},{(5.5,1)/b}, {(6.5,0)/c},{(5.5,-1)/d}}
       \node[pnode] (\name) at \pos {};
     \foreach \source/\dest in {a/d,b/d}
        \path[edge] (\source) -- (\dest);
\foreach \source/\dest in {a/b,c/b}
        \path[dedge] (\source) -- (\dest);
\foreach \source/\dest in {c/d}
        \path[tredge] (\source) -- (\dest);
    \foreach \pos/\name in {{(3,.8)/y}, {(3,-.8)/z}}
        \node[pnode] (\name) at \pos {};
     \foreach \source/\dest in {a/y,a/z}
        \path[edge] (\source) -- (\dest);
\foreach \source/\dest in {y/z}
        \path[quple] (\source) -- (\dest);
        \end{scope}
\end{tikzpicture}
\end{center}
\caption{Irreducible graphs with connectivity 1}\label{f:ircut}
\end{figure}

Since any cut vertex $v$ must have 3 edges to one component of $G-v$ and 2 edges to the other, there will be two irreducible graphs of connectivity 1 with six vertices as shown in figure \ref{f:ircut}. Any larger graph of connectivity 1 can be broken into two connected graphs which are each smaller than than the original graph. %will be reducible using at least one of the given reductions as will be shown 

In section \ref{s:cutred} we will show how to reduce any graph with a cut vertex adjacent to one of these minimal configurations, but until then we will be able to assume that any graph with vertex connectivity 1 is reducible if it has more than 6 vertices.
%and so is a graph with all vertices of degree 5 apart from one of degree 2 or 3.

\subsection{{Fundamental Graphs}}\label{s:irred}

The 4-regular graphs which are the line graphs of 3-regular graphs are used as a base family which cannot be reduced by the operations in \cite{ar:PR}, and we can generalise that idea as follows: 
an $(a,b)$-biregular graph is a bipartite graph without multiple edges in which all vertices in one part of the bipartition have degree $a$ and all other vertices have degree $b$.
Note that the line graph operation on a 3-regular graph can be viewed as first subdividing each edge to form a (2,3)-biregular graph $H$ and then deleting each vertex $v$ of degree 3 and replacing $v$ by a triangle through its neighbours in $H$.

A quintic graph $Q(B)$ can be formed from a (3,4)-biregular bipartite graph $B$ by creating one copy of $K_4$ for each of the vertices of degree 4 and, for each vertex $v$ of degree 3 in $B$, joining a set of three previously unused vertices from each $K_4$ corresponding to the neighbours of $v$. Such a graph will necessarily be quintic and have the triangle property; each edge in a $K_4$ is in triangles with the two other vertices from the $K_4$ and each edge between two copies of $K_4$ is in that one triangle. If $B$ had $n$ vertices of degree 4 then $Q(B)$ will necessarily have $4n$ vertices; however, since $B$ is biregular, there must also be $\frac{4n}{3}$ vertices of degree 3 in $B$, and hence $n$ must also be a multiple of 3. Therefore $Q(B)$ will have $12k$ vertices for some integer $k\geq1$ when $B$ has $3k$ vertices of degree 4 and $4k$ vertices of degree 3.

\begin{figure}[ht]
\begin{center}
\begin{tikzpicture}[scale=0.7, transform shape]
 \def\toa{4.5}
\def\tob{1}
\def\tta{1.4}
\def\ttb{4.45}
\def\tr{1.8}
%\node[math] () at (0.2,5) {$Q(K_{3,4}):$};
%\foreach \pos/\name in {{(7,5)/a1}, {(9,5)/b1}, {(10,6)/c1},{(10,3)/d1},{(3,4)/a2}, {(5,4)/b2}, {(6,5.5)/c2}, {(6,2)/d2}, {(6.5,.5)/a3}, {(8,0.5)/b3}, {(9,2)/c3},{(9,-1)/d3}}
      \foreach \pos/\name in {{(-0.866*\tr,-0.5*\tr)/a1}, {(0:0)/b1}, {(0,\tr)/c1},{(0.866*\tr,-0.5*\tr)/d1}}
        \node[enode] (\name) at \pos {};
          \foreach \pos/\name in {{(\toa-0.866*\tr,\tob-0.5*\tr)/a2}, {(\toa,\tob)/b2}, {(\toa,\tob+\tr)/c2},{(\toa+0.866*\tr,\tob-0.5*\tr)/d2}}
        \node[enode] (\name) at \pos {};
               \foreach \pos/\name in {{(\tta-0.866*\tr,\ttb-0.5*\tr)/a3}, {(\tta,\ttb)/b3}, {(\tta,\ttb+\tr)/c3},{(\tta+0.866*\tr,\ttb-0.5*\tr)/d3}}
        \node[enode] (\name) at \pos {};
    % First we draw the vertices
    % Connect vertices with edges and draw weights
    \foreach \source/\dest in {a1/b1, a1/c1,a1/d1,b1/c1,b1/d1,c1/d1}
        \path[edge] (\source) -- (\dest) [color=purple];
\foreach \source/\dest in {a2/b2, a2/c2,a2/d2,b2/c2,b2/d2,c2/d2}
        \path[edge] (\source) -- (\dest) [color=purple];
\foreach \source/\dest in {a3/b3, a3/c3,a3/d3,b3/c3,b3/d3,c3/d3}
        \path[edge] (\source) -- (\dest) [color=purple];
 \foreach \source/\dest in {a1/a2, a1/a3,a2/a3} \path[edge] (\source) -- (\dest) [dashed];
\foreach \source/\dest in {b1/b2, b1/b3,b2/b3} \path[edge] (\source) -- (\dest) [dashed];
\foreach \source/\dest in {c1/c2, c1/c3,c2/c3} \path[edge] (\source) -- (\dest) [dashed];
\foreach \source/\dest in {d1/d2, d1/d3,d2/d3} \path[edge] (\source) -- (\dest) [dashed];
\end{tikzpicture}
\end{center}
\caption{A graph formed from a $(3,4)$-biregular graph: $Q(K_{3,4})$}\label{f:k34}
\end{figure}
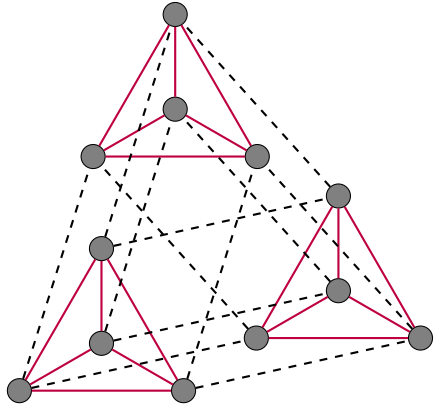

For instance, given $K_{3,4}$, we form the graph shown in figure \ref{f:k34}, where the $K_4$ edges are shown in purple and the triangle edges are dashed..
There are 18 different (3,4) biregular graphs with 6 vertices of degree 4 (created via nauty \cite{ar:nauty}) and so there are 18 non-isomorphic fundamental quintic graphs with 24 vertices formed in this way, and an infinite family of these graphs. However, they will be shown reducible in section \ref{s:final}.

Note that these graphs are made from edge-disjoint copies of $K_3$ and $K_4$ such that every vertex is in exactly one of each of these graphs. In general, all %2-connected 
irreducible quintic graphs with the triangle property will be shown to be constructed in this way from a small pool of configurations, which we shall call {\em atoms}, and $A_1:=K_3$ and $A_2:=K_4$ are the first atomic graphs, and the various ways they can be combined will be studied in section \ref{s:atom}.

\section{{General Reductions}}\label{s:redlist}

We will now introduce some basic reductions that we will be using to show that most 5-regular graphs with the triangle property can be formed recursively. 
%Suppose $G$ is such a graph. 
An eligible triangle was defined in \cite{ar:PR} as one that can be removed and the resulting graph keeps the triangle property, or if it had a triple edge. We will need  a different classification of triangles for our reductions.

Working in a specific graph $G$, let $m(T)$ be the number of edges in a triangle $T$ which %are not multiple edges and 
 are only in one triangle; note that since we are considering triangles as sets of vertices, multiple edges will not contribute extra triangles.
  If $m(T)=3$ then $T$ is eligible since no other triangles are affected by its removal, however, it is possible for a triangle with $m(T)<3$ to also be eligible in a quintic graph. We will refer to triangles with $m(T)=3$ as {\em aloof}.
 By observation \ref{o:four} any aloof triangles either have no multiple edges, one double edge or are the triangle with a quadruple edge from figure \ref{f:cutv}. These configurations are all atoms; $A_1$, $A_3$ and $A_4$, respectively.
 
 If $m(T)=2$ then we will call that triangle {\em unsafe} since removal of an edge from it will mean that there is now an edge which is not in a triangle, unless the edge removed was one of a multiple edge.
 When $m(T)=1$, if the edge of $T$ that is in only one triangle is $e$ then $G-e$ will still have the triangle property, but deletion of the other edges of $T$ (if they are not multiple) will mean that $e$ is no longer in a triangle. 
 %We will refer to such triangles as {\em fragile}.

\subsection{Z-reduction}

\begin{figure}[ht]
\begin{center}
\begin{tikzpicture}[scale=0.7, transform shape]
   \foreach \pos/\name in {{(0,0)/a}, {(0,2)/b}, {(2,2)/c},{(2,0)/d}, {(6,0)/e}, {(6,2)/f}}
        \node[enode] (\name) at \pos {};
   \foreach \pos/\name in {{(-.7,-.2)/a1}, {(-.2,-.7)/a2}, {(2.2,2.7)/c2},{(2.7,2.2)/c1}, {(-.2,2.7)/b2},{(-.7,2.2)/b1},{(-.5,2.5)/b3}, {(2.2,-.7)/d2},{(2.7,-.2)/d1},{(2.5,-.5)/d3},{(6-.7,-.2)/e1},{(6-.2,-.7)/e2},{(6.2,2.7)/f2},{(6.7,2.2)/f1},{(6-.2,2.7)/f3},{(6-.7,2.2)/f4},{(6-.5,2.5)/f5},{(6.2,-.7)/e5},{(6.7,-.2)/e4},{(6.5,-.5)/e3}}
        \node[nvertex] (\name) at \pos {};
    \foreach \source/\dest in {a/b,b/c,c/d,d/a,a/a1,a/a2,b/b1,b/b2,b/b3,c/c1,c/c2,d/d1,d/d2,d/d3}
        \path[edge] (\source) -- (\dest);
\foreach \source/\dest in {e/e1,e/e2,e/e3,e/e4,e/e5,f/f1,f/f2,f/f3,f/f4,f/f5}
        \path[edge] (\source) -- (\dest);
    \foreach \source/\dest in {a/c}
        \path[tedge] (\source) -- (\dest);
\node[math] () at (.5,1.3) {$H$};
\node[math] () at (.3,2.3) {$a$};
\node[math] () at (2.3,1.7) {$b$};
\node[math] () at (-.3,.3) {$c$};
\node[math] () at (1.7,-.3) {$d$};
\node[math,below] () at (f.south) {$\{a,b\}$};
\node[math,above] () at (e.north) {$\{c,d\}$};
\node[math] () at (4,1) {$\tikzLongrightarrow$};
\end{tikzpicture}
\end{center}
\caption{The Z-reduction $abcd$}\label{f:zredi}
\end{figure}
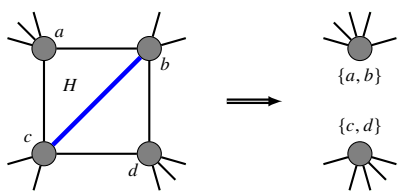

Let $G$ be a quintic graph with the triangle property. By corollary \ref{c:triv} we know there will likely be many diamonds, unless $G$ is an extreme case such as those introduced in subsection \ref{s:irred}.
 If a diamond exists which is not part of a $K_4$ in $G$ (and has no double edges), we can remove its five edges and identify together two pairs of vertices of degree two and three which were previously in the diamond as shown in figure \ref{f:zredi}. That is; we can either identify $a$ with $b$ and $c$ with $d$, or $a$ with $c$ and $b$ with $d$; we shall refer to the vertices in the order of the underlying Z shaped path, so the Z-reduction is either $abcd$ or $acbd$.

This operation will preserve the triangle property if none of the edges of the diamond were multiple or used as parts of triangles outside of this configuration.  
For the Z-reduction using $abcd$ on the diamond $H$, the following properties will mean that the Z-reduction does not produce a 5-regular graph with the triangle property:
\begin{itemize}
\item[Z1:]
If $ab$ or $cd$ is a multiple edge then a loop is formed, which cannot be in a triangle. If any of the other edges are multiple then we may be able to use the X-reduction which will be introduced in section \ref{ss:xred}.
\item[Z2:]
If $ab$ or $cd$ are part of an unsafe triangle $T$ outside of $H$ then, on reduction, $T$ will collapse into a multiple edge and that edge will not be in a triangle if both of the other edges of $T$ were only in that triangle, i.e. $T$ is unsafe.%, or $m(T)=2$.
\item[Z3:]
If any of $ac$, $bc$ or $bd$ are part of a triangle $T$ outside of $H$ with $m(T)\geq1$ then the reduced graph will not have the triangle property.
\end{itemize}

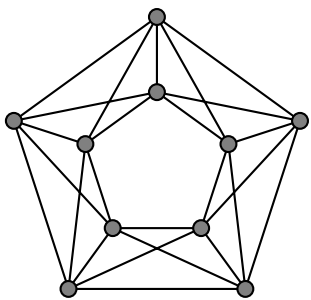
\begin{figure}[ht]
\begin{center}
\begin{tikzpicture}[thick,rotate=18,scale=0.5]
    \draw \foreach \x in {0,72,...,288}
    {
    %polar coordinates (theta:r)
          (\x:2) node {} -- (\x:4)
	   (\x:4) node {}  -- (\x+72:4)
  	   (\x:2) -- (\x+72:2)
  	   (\x:4) -- (\x+72:2)
  	   (\x:2) -- (\x+72:4)
        };
     \foreach \x in {0,72,...,288} {\node[snode]  at (\x:2) {};\node[snode]  at (\x:4) {};}
   %   \foreach \x in {0,72,...,288} 
\end{tikzpicture}
\end{center}
\caption{A graph with no Z-reductions possible}\label{f:k4cycle}
\end{figure}

%In particular, for the graphs in section \ref{s:intro} such as $L_m$ this reduction will always work since they are simple and the only triangles are those in the diamonds. 
%In general, we know there is at least one diamond subgraph at every vertex without multiple edges and two ways to apply the Z-reduction at each diamond there are many possible ways to apply this reduction to a graph. 
Most graphs do have a diamond subgraph which can be used to reduce it using the Z-reduction, but
there do exist some graphs for which it does not give a 5-regular graph with the triangle property, such as the one in figure \ref{f:k4cycle}.

Most diamonds in the graph in figure \ref{f:k4cycle} are part of a $K_4$, but there are some edges (the radial spokes) which are not in a $K_4$ but are part of four triangles and hence four diamonds also. The resulting graph from a Z-reduction centred on those edges does not have the triangle property.

\subsection{\texorpdfstring {$K_4$}{K4} subgraphs or multiple edges in a diamond}\label{ss:xred}

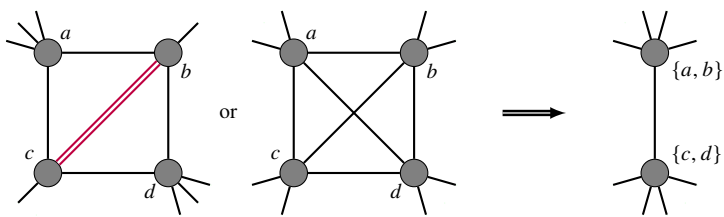
\begin{figure}[ht]
\begin{center}
\begin{tikzpicture}[scale=0.8, transform shape]
   \foreach \pos/\name in {{(0,0)/a}, {(0,2)/b}, {(2,2)/c},{(2,0)/d}}
        \node[enode] (\name) at \pos {};
   \foreach \pos/\name in {{(-.5,-.5)/a1}, {(-.2,-.7)/a2}, {(2.5,2.5)/c2}, {(-.2,2.7)/b2},{(-.7,2.2)/b1},{(-.5,2.5)/b3}, {(2.2,-.7)/d2},{(2.7,-.2)/d1},{(2.5,-.5)/d3}}
        \node[nvertex] (\name) at \pos {};
    % First we draw the vertices
    % Connect vertices with edges and draw weights
    \foreach \source/\dest in {b/c,a/d,b/a,c/d,a/a1,b/b1,b/b2,b/b3,c/c2,d/d1,d/d2,d/d3}
        \path[edge] (\source) -- (\dest);
\foreach \source/\dest in {c/a}
        \path[dedge] (\source) -- (\dest);
\node[math] () at (3,1) {or};
\node[math] () at (.3,2.3) {$a$};
\node[math] () at (2.3,1.7) {$b$};
\node[math] () at (-.3,.3) {$c$};
\node[math] () at (1.7,-.3) {$d$};
\end{tikzpicture}
\begin{tikzpicture}[scale=0.8, transform shape]
   \foreach \pos/\name in {{(0,0)/a}, {(0,2)/b}, {(2,2)/c},{(2,0)/d}, {(6,0)/e}, {(6,2)/f}}
        \node[enode] (\name) at \pos {};
   \foreach \pos/\name in {{(-.7,-.2)/a1}, {(-.2,-.7)/a2}, {(2.2,2.7)/c2},{(2.7,2.2)/c1}, {(-.2,2.7)/b2},{(-.7,2.2)/b1},{(-.5,2.5)/b3}, {(2.2,-.7)/d2},{(2.7,-.2)/d1},{(2.5,-.5)/d3},{(6-.7,-.2)/e1},{(6-.2,-.7)/e2},{(6.2,2.7)/f2},{(6.7,2.2)/f1},{(6-.2,2.7)/f3},{(6-.7,2.2)/f4},{(6-.5,2.5)/f5},{(6.2,-.7)/e5},{(6.7,-.2)/e4},{(6.5,-.5)/e3}}
        \node[nvertex] (\name) at \pos {};
    % First we draw the vertices
    % Connect vertices with edges and draw weights
    \foreach \source/\dest in {a/a1,a/a2,b/b1,b/b2,c/c1,c/c2,d/d1,d/d2}
        \path[edge] (\source) -- (\dest);
\foreach \source/\dest in {e/f,e/e1,e/e2,e/e4,e/e5,f/f1,f/f2,f/f3,f/f4}
        \path[edge] (\source) -- (\dest);
    \foreach \source/\dest in {a/b,b/c,c/d,d/a,a/c,b/d}
        \path[edge] (\source) -- (\dest);
\node[math] () at (.3,2.3) {$a$};
\node[math] () at (2.3,1.7) {$b$};
\node[math] () at (-.3,.3) {$c$};
\node[math] () at (1.7,-.3) {$d$};
\node[math] () at (4,1) {$\tikzLongrightarrow$};
\node[math,right] () at (6.15,1.7) {$\{a,b\}$};
\node[math,right] () at (6.15,.3) {$\{c,d\}$};
\end{tikzpicture}
\end{center}
\caption{Multiple edge Z-reduction and X-reduction}\label{f:xred}
\end{figure}

It is possible to define a similar reduction to the Z-reduction for a $K_4$ subgraph as shown in figure \ref{f:xred}. This time, due to the symmetry of $K_4$ there can be up to three different pairs of vertices which can be identified, but it is now necessary for there to be at least one vertex outside the subgraph which is a neighbour of two vertices of the $K_4$ so that the edge between the two vertices in the reduced graph is part of a triangle. As shown, the Z-reduction also gives the same reduced graph when there is a multiple edge which is not contracted, and, more generally, if there was more than one such multiple edge then we can just form more multiple edges in the reduced graph. 

The following are the situations in which the multiple edge Z-reduction or X-reduction contracting $ab$ and $cd$ will not give a 5-regular graph with the triangle property, which are similar to those for the Z-reduction:
\begin{itemize}
\item[X1:]
If $ab$ or $cd$ is a multiple edge then a loop is again formed. 
\item[X2:]
If $ab$ or $cd$ are part of an unsafe triangle outside of the $K_4$.
\item[X3:]
If none of $ac$, $bc$ and $bd$ (and $ad$ for $K_4$) are part of a triangle $T$ outside of the configuration then the edge between the identified vertices in the reduced graph will not be in a triangle.
% with $m(T)=2$ then $T$ will collapse into a multiple edge and that edge will, as in Z2, not be in a triangle in the reduced graph.
% if both of the other edges of $T$ were only in that triangle.
\end{itemize}

%Graph1: $K_5$ -> disappears
%Graph2: $K_4$ -> cases
%Graph3: $W_5$+ -> double edges
%Graph4: 4-edge + any multiple - > disappears
%Graph5: 3-edge
%Graph6: Z+

\subsection{Large Complete subgraphs}

Suppose $G$ is a connected quintic graph with the triangle property, and its clique number is $\omega(G)$.
If $\omega(G)\geq 6$ then $G$ must be $K_6$ and no multiple edges are possible. We are able to prove a similar result to Lemma 2 in \cite{ar:PR} which will be useful in this section:

\begin{lem}\label{l:eout} 
An induced subgraph $H$ of $G$ with between one and three vertices of degree 4 and the others of degree 5 must have all edges from $H$ incident with a single vertex in $G$, which will be a cut-vertex.
\end{lem}

\begin{proof}
We need to consider the set $S$ of edges from $H$ to $G\backslash V(H)$. Each must be in a triangle, and so must have a vertex in common with another edge in $S$, but all vertices in $H$ are incident with at most one edge in $S$, so edges in $S$ share a vertex outside of $V(H)$. Since $|S|\leq 3$ and we need two edges per vertex to be in the triangle we cannot have more than one vertex joining $S$ to $G$, so it must be a cut-vertex. Additionally, $|S|\geq 2$ because a single edge can't be in a triangle. 
\end{proof}

\begin{thm}\label{t:cn5}
All quintic graphs with the triangle property which have clique number 5 are reducible
\end{thm}

\begin{proof}
Suppose $G$ is a quintic graph with the triangle property and $\omega(G)=5$, and let $H$ be an induced subgraph of $G$ containing five mutually adjacent vertices. If $|E(H)|>10$ then there are at most 3 edges from $H$ to $G$ and thus, by 
lemma \ref{l:eout}, $H$ must be adjacent to a cut-vertex. Since $H$ has five vertices we can replace it by a smaller graph using a reduction from figure \ref{f:cutv}. 

If $|E(H)|=10$ then each vertex in $H$ must be joined to one vertex outside of $H$. Since all of the edges from $H$ to these vertices must be in triangles there must either be one vertex joined to all of $H$ (which would make $\omega(G)>5$, contrary to our supposition) or there are two vertices, which act similarly to those in lemma \ref{l:eout}. One vertex $u$ is joined by two edges to $H$ and the other $v$ is joined by three edges.
A double edge at $v$ into $H$ is not possible since that vertex could then not be joined to all four other vertices of $H \cong K_5$.

It is possible to use the X-reduction for this subgraph; let $u_1$ and $u_2$ be the neighbours of $u$ and $v_1$ and $v_2$ two of the neighbours of $v$. We contract $u_1v_1$ and $u_2v_2$ and the vertices resulting from the contraction will form triangles with both $u$ and $v$ and the fifth vertex of $H$ will have double edges to both contracted vertices and a single edge to $u$, all of which are in triangles and each vertex is degree 5 as required. Note that none of the properties X1, X2 or X3 can hold 
since $H$ has no multiple edges and none of the triangles outside of the $K_4$ in $H$ are unsafe.
\end{proof}

%\subsection{Clique Number Four}\label{s:cn4}

We can, from now on, suppose that $G$ is a quintic graph with the triangle property and that $\omega(G)\leq 4$; for the remainder of this section we will suppose $\omega(G)=4$, so $G$ has at least one subgraph isomorphic to $K_4$.

\begin{thm}\label{t:cn4a}
All quintic graphs with the triangle property which have clique number 4 and a vertex adjacent to three vertices of a 4-clique are reducible.
\end{thm}

\begin{proof}
Let $H$ be a subgraph of $G$ isomorphic to $K_5$ with one edge removed as shown %induced by $\{u,v,w_1,w_2,w_3\}$
 on the left in in figure \ref{f:k5me}.
 There will be at most seven edges joining $H$ to $G\backslash V(H)$, but
as in theorem \ref{t:cn5} we can use lemma \ref{l:eout} to simplify the situation to either there being one double edge in $H$ or $H$ being simple. 

We can suppose there is not a double edge at $w_1$, say, as there are three symmetrical vertices 
$\{w_1,w_2,w_3\}$. By the symmetry between $u$ and $v$, we can also suppose there is no double edge at $u$ so it must be adjacent to two vertices outside of $H$ since if $uv$ was an edge $\omega(G)=5>4$.
Thus there are only two possible double edges in $H$, either $w_2w_3$ or $vw_2$, without loss of generality.
%as shown in figure \ref{f:k5me} 

\begin{figure}[ht]
\begin{center}
\begin{tikzpicture}[scale=0.8, transform shape]
   \foreach \pos/\name in {{(0,0)/a}, {(4,0)/b}, {(2,3)/c},{(2,.7)/d}, {(2,1.4)/e}}
        \node[enode] (\name) at \pos {};

%        \node[nvertex] (\name) at \pos {};
    % First we draw the vertices
    % Connect vertices with edges and draw weights
    \foreach \source/\dest in {a/b,b/c,c/a,a/d,a/e,c/e,b/d,b/e,d/e}
        \path[edge] (\source) -- (\dest);
%\foreach \source/\dest in {e/f,e/e1,e/e2,e/e4,e/e5,f/f1,f/f2,f/f3,f/f4}
   %     \path[edge] (\source) -- (\dest);
%\node[math] () at (0,1.4) {$H$};
\node[math] () at (2.4,3) {$u$};
\node[math] () at (2,0.2) {$v$};
\node[math] () at (-.5,0) {$w_1$};
\node[math] () at (1.5,1.5) {$w_2$};
\node[math] () at (4.5,0) {$w_3$};
\node[math] () at (5,1.5) {$\tikzLongrightarrow$};

   \foreach \pos/\name in {{(6,0)/a}, {(10,0)/b}, {(8,3)/c}}
        \node[enode] (\name) at \pos {};
 \foreach \source/\dest in {a/c,b/c}
        \path[edge] (\source) -- (\dest);
 \foreach \source/\dest in {a/b}
        \path[dedge] (\source) -- (\dest);
\node[math] () at (8.8,3) {$\{u,w_3\}$};
\node[math] () at (10.5,0) {$v$};
\node[math] () at (6,-.5) {$\{w_1,w_2\}$};
\end{tikzpicture}
\end{center}
\caption{A $K_4$ with a neighbour adjacent to three of its vertices, reduced to $A_3$}\label{f:k5me}
\end{figure}
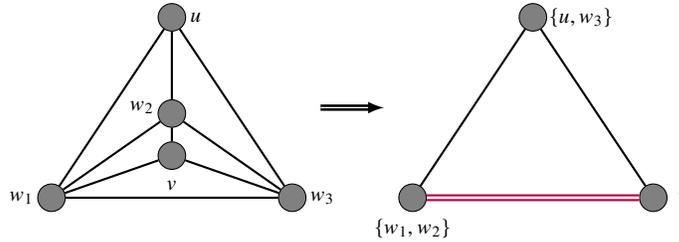

We choose to perform an X-reduction by contracting $uw_3$ and $w_1w_2$ and a quintic graph $G'$ will be formed. Since neither $u$ nor $w_1$ is part of a double edge, property X1 does not apply. X3 does not apply since there is a triangle containing $v$, $w_2$ and $w_3$. 

For property X2 to apply there would need to be a vertex outside of $H$ which was adjacent to both $u$ and $w_3$ or both $w_1$ and $w_2$ which gave rise to an unsafe triangle. However, there is another similar X-reduction possible, contracting $w_1w_3$ and $uw_2$ instead to form $G''$, giving the same basic structure as in figure \ref{f:k5me}, and this also can only have property X2, but this time if there is an unsafe triangle using the edges $uw_2$ or $w_1w_3$.

Firstly suppose there is a double edge in $H$. By the symmetry discussion earlier the double edge is from $w_2$ and that would be the fifth known edge from $w_2$.  Thus $w_2$ couldn't have a common neighbour with another vertex and we only need to worry for property X2 about whether $uw_3$ or $w_1w_3$ are part of an unsafe triangle. Since both of these edges involve $w_3$, we can assume that the double edge was not $w_2w_3$, so must be $vw_2$. 
However, $w_3$ has only one unknown neighbour, say $y$, so we can choose to use either $G'$ or $G''$ unless $y$ is adjacent to both $u$ and $w_1$, but in such a case the triangles formed will have $m(T)<2$ as $\{u,w_1,w_3,y\}$ induce a $K_4$ and so all of these edges are in more than one triangle.

If there is not a double edge in $H$ then, as in the double edge case, we are concerned only if neither $G'$ or $G''$ give a valid reduction, which means that we would require two unsafe triangles adjacent in the 4-cycle $uw_3w_1w_2$. Since all $w_j$ vertices already have four known neighbours, any triangle involving them and a new vertex could never be unsafe, and for $u$ the same $K_4$ as before will exist. 

Thus $H$ is reducible in every case.
\end{proof}

From henceforth in this section, we can assume that, in $G$, any subgraph $H$ isomorphic to $K_4$ has no vertex outside of $V(H)$ adjacent to more than two of its vertices. We will break the cases down by a parameter $s_H$ which is defined as the number of vertices in $G\backslash H$ adjacent to two vertices of $H$.

There are four vertices in $H$ of degree at least 3, so, as $G$ is quintic, there are at most eight edges available to join $H$ to the set $S$ of vertices adjacent to $H$, so $s_H\leq \frac{8}{2} = 4$. Moreover, if there are $k$ multiple edges in $H$ then each extra edge reduces the possible number of vertices in $S$ and so $s_H \leq 4-k$.
If $k\geq 3$ then $s_H \leq 4-3 \leq 1$ and if, additionally, $s_H=1$ then there is a unique graph with a cut vertex which can be reduced as in figure \ref{f:cutv}. Thus we can assume that there are at most two double edges or one triple edge until case iii) when $s_H=0$.

\begin{itemize}
\item{Case i) $3\leq s_H \leq4$}\\
$H$ must be simple to have four vertices outside of $H$ in triangles with edges of $H$ and so property X1 cannot be satisfied when $s_H=4$. Since $\binom{4}{2}=3$ we have three different choices of pairs of edges in $H$ to contract in the X-reduction so the only way that property X2 or X3 can be satisfied for all possible X-reductions is if there is a triangle in $H$ which has all three pairs of its vertices adjacent to vertices in $S$. 
%This implies that $s_H=3$ and 
We therefore have the configuration shown in figure \ref{f:sh3}, and, moreover, $t$ is not joined to $u$, $v$ or $w$ since that would give a vertex adjacent to three vertices of a $K_4$.
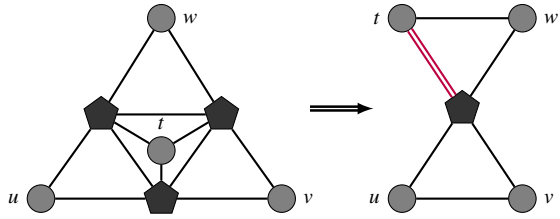
\begin{figure}[ht]
\begin{center}
\begin{tikzpicture}[scale=0.8, transform shape]
   \foreach \pos/\name in {{(0,0)/a}, {(4,0)/b}, {(2,3)/c},{(2,.8)/d}}
        \node[enode] (\name) at \pos {};
  \foreach \pos/\name in {{(3,1.4)/e}, {(1,1.4)/f}, {(2,0)/g}}
        \node[pnode] (\name) at \pos {};

%        \node[nvertex] (\name) at \pos {};
    % First we draw the vertices
    % Connect vertices with edges and draw weights
    \foreach \source/\dest in {a/g,b/g,b/e,e/c,c/f,f/a,f/g,g/e,e/f,d/f,d/e,d/g}
        \path[edge] (\source) -- (\dest);
\node[math,left] () at (a.west) {$u$};
\node[math,right] () at (b.east) {$v$};
\node[math,right] () at (c.east) {$w$};
\node[math,above] () at (d.north) {$t$};
%\foreach \source/\dest in {e/f,e/e1,e/e2,e/e4,e/e5,f/f1,f/f2,f/f3,f/f4}
   %     \path[edge] (\source) -- (\dest);
%\node[math] () at (0,1.4) {$H$};
\node[math] () at (5,1.5) {$\tikzLongrightarrow$};

   \foreach \pos/\name in {{(6,0)/a}, {(8,0)/b}, {(6,3)/d}, {(8,3)/e}}
        \node[enode] (\name) at \pos {};
   \foreach \pos/\name in {{(7,1.5)/c}}
        \node[pnode] (\name) at \pos {};
 \foreach \source/\dest in {a/c,b/c,a/b,c/e,d/e}
        \path[edge] (\source) -- (\dest);
 \foreach \source/\dest in {c/d}
        \path[dedge] (\source) -- (\dest);
\node[math,left] () at (a.west) {$u$};
\node[math,right] () at (b.east) {$v$};
\node[math,right] () at (e.east) {$w$};
\node[math,left] () at (d.west) {$t$};
\end{tikzpicture}
\end{center}
\caption{No X-reduction possible when $\omega(G)=4$}\label{f:sh3}
\end{figure}
However, it is possible to use  Z-reduction on any of the outer diamonds to give the resulting configuration which is quintic and has the triangle property since none of the remaining edges from the named vertices were parts of triangles with edges that were deleted. There cannot be a multiple edge in this configuration since all edges have at least one vertex of degree 5. If $s_H=3$ and there is a double edge, then the X-reduction can be used.

\item{Case ii) $1\leq s_H \leq2$}\\
If $H$ is simple then property X1 cannot hold, and, similarly to case i), because $s_H<3$ we can choose a pair of non incident edges in $H$ which are not part of a triangle outside $H$, so X2 cannot hold for them. A triangle guaranteed by $s_H>0$ will mean X3 cannot hold.

Since there are fewer than 3 multiple edges in $H$ we can choose a pair of edges in $H$ to use for an X-reduction without satisfying property X1. A triangle from $s_H>0$ will mean that one of X3 and X2 can only hold if in $H$ there is a triangle which has either two double edges and one unsafe triangle outside $H$ or one double edge and two unsafe triangles outside $H$. 

Both of these situations can be reduced to a quintic graph by deleting some pentagonal vertices and identifying a vertex of degree 3 with one of degree 2 as shown in figure \ref{f:ident23}; they are guaranteed to not be adjacent since there is no vertex adjacent to three vertices of a $K_4$.
The resulting graphs will have the triangle property since $u_1$ and $u_2$ cannot be joined by an edge as the triangles joining them to $H$ must be unsafe.

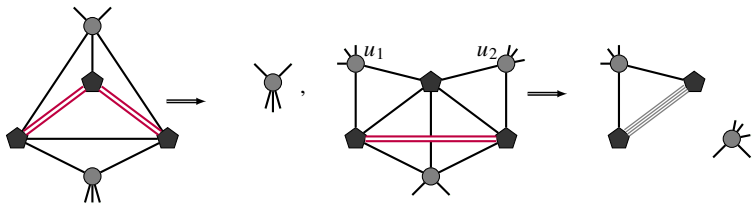
\begin{figure}[ht]
\begin{center}
\begin{tikzpicture}[scale=0.5, transform shape]
\begin{scope}
   \foreach \pos/\name in {{(2,3)/d},{(2,-1)/e}}
        \node[enode] (\name) at \pos {};
          \foreach \pos/\name in {{(0,0)/a}, {(4,0)/b}, {(2,1.5)/c}}
        \node[pnode] (\name) at \pos {};
        \foreach \pos/\name in {{(2.5,3.5)/d1}, {(1.5,3.5)/d2}, {(2.2,-1.7)/e2},{(1.8,-1.7)/e1},{(2,-1.7)/e3}}
        \node[nvertex] (\name) at \pos {};
    \foreach \source/\dest in {a/b,e/a,e/b,d/a,d/b,d/c,d/d1,d/d2,e/e1,e/e2,e/e3}
        \path[edge] (\source) -- (\dest);
\foreach \source/\dest in {b/c,c/a}
        \path[dedge] (\source) -- (\dest);
        \node[math,rotate=0] () at (4.5,1) {$\tikzLongrightarrow$};
\end{scope}
\begin{scope}[shift={(9,0)}]
   \foreach \pos/\name in {{(0,2)/d},{(2,-1)/e},{(4,2)/f}}
        \node[enode] (\name) at \pos {};
          \foreach \pos/\name in {{(0,0)/a}, {(4,0)/b}, {(2,1.5)/c}}
        \node[pnode] (\name) at \pos {};
        \foreach \pos/\name in {{(-.5,2)/d1}, {(-.3,2.4)/d2}, {(0,2.5)/d3},{(2.5,-1.5)/e1},{(1.5,-1.5)/e2}, {(4.1,2.5)/f3}, {(4.5,2.1)/f2},{(4.3,2.4)/f1}}
        \node[nvertex] (\name) at \pos {};
    \foreach \source/\dest in {a/c,b/c,c/e,e/a,e/b,d/a,f/b,f/c,d/c,d/d1,d/d2,d/d3,e/e1,e/e2,f/f1,f/f2,f/f3}
        \path[edge] (\source) -- (\dest);
\foreach \source/\dest in {b/a}
        \path[dedge] (\source) -- (\dest);
        \node[math] () at (0.5,2.3) {{\LARGE $u_1$}};
        \node[math] () at (3.5,2.3) {{\LARGE $u_2$}};
        \node[math,rotate=0] () at (5.1,1.2) {$\tikzLongrightarrow$};
\end{scope}
\begin{scope}[shift={(4.8,0)}]
   \foreach \pos/\name in {{(2,1.5)/c}}
        \node[enode] (\name) at \pos {};
        \foreach \pos/\name in {{(2.5,2)/d1}, {(1.5,2)/d2}, {(2.2,0.8)/e2},{(1.8,0.8)/e1},{(2,0.7)/e3}}
        \node[nvertex] (\name) at \pos {};
    \foreach \source/\dest in {c/d1,c/d2,c/e1,c/e2,c/e3}
        \path[edge] (\source) -- (\dest);
        \node[math,rotate=0] () at (2.8,1.2) {{\LARGE ,}};
        \end{scope}
\begin{scope}[shift={(16,0)}]
   \foreach \pos/\name in {{(0,2)/d},{(3,0)/e}}
        \node[enode] (\name) at \pos {};
          \foreach \pos/\name in {{(0,0)/a}, {(2,1.5)/c}}
        \node[pnode] (\name) at \pos {};
        \foreach \pos/\name in {{(-.5,2)/d1}, {(-.3,2.4)/d2}, {(0,2.5)/d3},{(3.5,-.5)/e1},{(2.5,-.5)/e2}, {(3.1,0.5)/f3}, {(3.5,0.1)/f2},{(3.3,0.4)/f1}}
        \node[nvertex] (\name) at \pos {};
    \foreach \source/\dest in {d/a,d/c,d/d1,d/d2,d/d3,e/e1,e/e2,e/f1,e/f2,e/f3}
        \path[edge] (\source) -- (\dest);
\foreach \source/\dest in {c/a}
        \path[quple] (\source) -- (\dest);
\end{scope}
\end{tikzpicture}
\end{center}
\caption{Identifying two vertices to reduce configurations}\label{f:ident23}
\end{figure}
%The edges in $H$ which are in the triangles between $H$ and $G\backslash H$ cannot be used for X-reductions, but there is at least one copy of $2K_2$ in $H$ available for an X-reduction.

\item{Case iii) $s_H=0$}\\
 Since there are no vertices in $G\backslash H$ adjacent to two vertices of $H$ then $S$ contains eight different vertices in four triangles with $H$ if there are no multiple edges in $H$. In the simple case any application of an X-reduction will leave an edge not in a triangle as in property X3; this structure was the basis of the fundamental graphs introduced in subsection \ref{s:irred}.
 
 If there is a multiple edge in $H$, then it is possible that there are no vertices at all in $G\backslash H$, in which case $G$ is one of the two graphs in figure \ref{f:fourv}, which correspond to adding either $C_4$ or $2C_2$ to $H$ as multiple edges.
 Similarly, we can add other 2-regular multigraphs as multiple edges to $H$; in this case we could have three double edges in a triangle in $H$, giving the left subgraph with a cut-vertex in figure \ref{f:cutv} or a triple edge between two vertices of $H$, as shown in figure \ref{f:trik4}, which will refer to as atom $A_5$.

\begin{figure}[ht]
\begin{center}
\begin{tikzpicture}[scale=0.5, transform shape]
   \foreach \pos/\name in {{(2,0)/b},{(2,4)/d}}
        \node[pnode] (\name) at \pos {};
           \foreach \pos/\name in {{(0,2)/a}, {(4,2)/c}}
        \node[enode] (\name) at \pos {};
   \foreach \source/\dest in {a/c,a/b,c/b,d/c,a/d}
        \path[edge] (\source) -- (\dest);
    \foreach \source/\dest in {b/d}
        \path[tredge] (\source) -- (\dest);
          \node[math] () at (-1,2) {{\LARGE $A_5$:}};
 \end{tikzpicture}
\end{center}
\caption{Atomic $K_4$ with an triple edge}\label{f:trik4}
\end{figure}
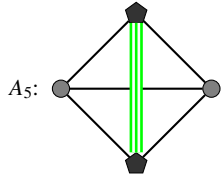

 However, if there is a vertex incident with a multiple edge in $H$ joined to a vertex $z$ outside of $H$ then the edge to $z$ must be in a triangle and so $z$ must be adjacent to another vertex of $H$ to accomplish this, contradicting $s_H=0$, so the graphs above are the only ones with $s_H=0$.

\end{itemize}

Thus the only quintic graphs with the triangle property with clique number 4 that we cannot yet reduce are those in which
all $K_4$ subgraphs have all vertices in their open neighbourhood only adjacent to one vertex of the $K_4$ and there is either one triple edge or no multiple edges in the $K_4$. Note that both of these atoms share a similar property regarding the triangles which share a vertex with their vertices (we shall call these triangles {\em pendant}).

\begin{lem}\label{l:a5red}
Atom $A_5$ is reducible and $A_2$ is reducible unless all of its pendant triangles are aloof.
\end{lem}

\begin{proof}
Since $s_H=0$ for both of these configurations, we know that all of the pendant triangles $T_j$ necessarily have $m(T_j)\geq 2$. There are two reductions for such triangles as shown in figure \ref{f:ptri}. In general, we delete all pentagonal vertices and then if a triangle is aloof, as on the left of the top $A_5$ in figure \ref{f:ptri}, we can remove all edges of the pendant triangle and rejoin them to two new vertices, which are then joined by enough edges to make them degree 5. If the triangle is not aloof, but unsafe as on the right of the top $A_5$, we can add another edge between the two vertices of the adjacent triangle.

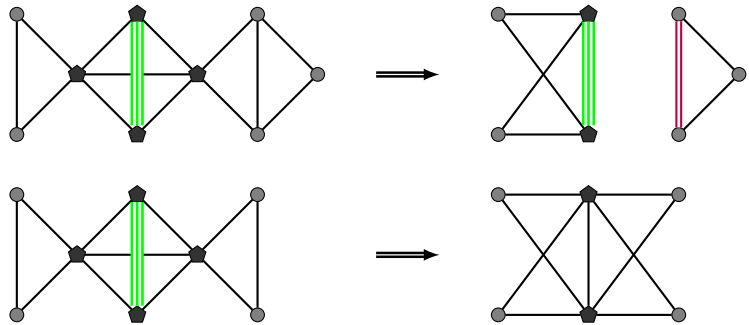
\begin{figure}[ht]
\begin{center}
\begin{tikzpicture}[scale=0.8, transform shape]
 \begin{scope}[shift={(0,0)},scale=0.5]
   \foreach \pos/\name in {{(2,0)/b},{(2,4)/d},{(0,2)/a}, {(4,2)/c}}
        \node[pnode] (\name) at \pos {};
           \foreach \pos/\name in {{(-2,0)/e}, {(-2,4)/f},{(6,0)/g}, {(6,4)/h},{(8,2)/i}}
        \node[enode] (\name) at \pos {};
   \foreach \source/\dest in {a/c,a/b,c/b,d/c,a/d,e/a,e/f,f/a,c/g,c/h,g/h,g/i,h/i}
        \path[edge] (\source) -- (\dest);
    \foreach \source/\dest in {b/d}
        \path[tredge] (\source) -- (\dest);
           \end{scope}
         \node[math,rotate=0] () at (5.5,1) {$\tikzLongrightarrow$};
         \begin{scope}[shift={(7.5,0)},scale=0.5]
            \foreach \pos/\name in {{(2,0)/b},{(2,4)/d}}
        \node[pnode] (\name) at \pos {};
           \foreach \pos/\name in {{(-1,0)/e}, {(-1,4)/f},{(5,0)/g}, {(5,4)/h},{(7,2)/i}}
        \node[enode] (\name) at \pos {};
   \foreach \source/\dest in {b/f,e/b,e/d,f/d,g/i,h/i}
        \path[edge] (\source) -- (\dest);
    \foreach \source/\dest in {b/d}
        \path[tredge] (\source) -- (\dest);
          \foreach \source/\dest in {g/h}
        \path[dedge] (\source) -- (\dest);
        \end{scope}
         \begin{scope}[shift={(0,-3)},scale=0.5]
   \foreach \pos/\name in {{(2,0)/b},{(2,4)/d},{(0,2)/a}, {(4,2)/c}}
        \node[pnode] (\name) at \pos {};
           \foreach \pos/\name in {{(-2,0)/e}, {(-2,4)/f},{(6,0)/g}, {(6,4)/h}}
        \node[enode] (\name) at \pos {};
   \foreach \source/\dest in {a/c,a/b,c/b,d/c,a/d,e/a,e/f,f/a,c/g,c/h,g/h}
        \path[edge] (\source) -- (\dest);
    \foreach \source/\dest in {b/d}
        \path[tredge] (\source) -- (\dest);
           \end{scope}
         \node[math,rotate=0] () at (5.5,-2) {$\tikzLongrightarrow$};
         \begin{scope}[shift={(7.5,-3)},scale=0.5]
            \foreach \pos/\name in {{(2,0)/b},{(2,4)/d}}
        \node[pnode] (\name) at \pos {};
           \foreach \pos/\name in {{(-1,0)/e}, {(-1,4)/f},{(5,0)/g}, {(5,4)/h}}
        \node[enode] (\name) at \pos {};
   \foreach \source/\dest in {b/f,e/b,e/d,f/d,d/b,g/d,h/d,g/b,h/b}
        \path[edge] (\source) -- (\dest);
        \end{scope}
\end{tikzpicture}
\end{center}
\caption{Pendant triangle reductions for $A_5$}\label{f:ptri}
\end{figure}

If, attached to $A_5$, there are two aloof triangles (as shown in the bottom reduction in figure \ref{f:ptri}) we will only need to add one edge between the two new vertices, and thus any combination of pendant triangles attached to $A_5$ allows a reduction. For $A_2$ the similar reductions will work if there are less than three aloof triangles in the same way since only two new vertices will need to be used as for $A_5$ and four pentagonal vertices were deleted.

If there are three aloof triangles attached to $A_2$ then it is possible to delete all four vertices of $A_2$, double the remaining edge of the unsafe triangle and then join the six vertices from the aloof triangles into two new aloof triangles, giving all vertices degree 5 and still having the triangle property. Thus the only currently irreducible $K_4$ configurations appear in graphs with all aloof pendant triangles such as those introduced in subsection \ref{s:irred}.
\end{proof}

\section{Clique Number Three}\label{s:reduc}

In this section we can suppose that there are no $K_4$ subgraphs in our quintic graph $G$ which has every edge in a triangle. We will mainly be using the Z-reduction, but also the X-reduction when the diamond configuration contains multiple edges. 
We will be able to use some special reductions in cases where the configurations have a number of vertices of degree 2, as well as, perhaps, some multiple edges. This generalises the reduction shown in the bottom of figure \ref{f:ptri}.

\begin{lem}\label{l:2s}
Given a quintic graph $G$ with the triangle property, suppose we have a configuration $H$ with only vertices of degree 5 and 2.
If the number of vertices of degree 2 is congruent to $c\mod{3}$ then $H$ can be reduced if there are $c$ multiple edges in $G$.
\end{lem}

\begin{proof}
Given such a configuration $H$ we proceed by deleting all the vertices of degree 5 and any edges between vertices of degree 2
(note that vertices of degree 2 could only be adjacent to each other in aloof triangles so no external triangles are affected).
 We need to add two edges to each of the vertices that were degree 2, and we can do that by using aloof triangles, which will preserve the triangle property and ensure the remaining graph is quintic if $c=0$. 

\begin{figure}[ht]
\begin{center}
\begin{tikzpicture}[scale=0.6, transform shape]
   \foreach \pos/\name in {{(-1,0)/c},{(-1,2)/d}, {(1,0)/e}, {(1,2)/f}, {(2,0)/g}, {(2,2)/h}}
        \node[enode] (\name) at \pos {};
  \foreach \pos/\name in {{(0,0)/a}, {(0,2)/b}}
        \node[pnode] (\name) at \pos {};
    \foreach \source/\dest in {a/b, a/c,a/d,a/e,a/f, b/c,b/d,b/e,b/f}
        \path[edge] (\source) -- (\dest);
           \foreach \source/\dest in {g/h}
        \path[dedge] (\source) -- (\dest);
        \node[math,right] () at (2.2,1) {\Large $e$};
     
      \node[math] () at (4,1) {$\tikzLongrightarrow$};
         \begin{scope}[shift={(8,0)}]
      \foreach \pos/\name in {{(-1,0)/c},{(-1,2)/d}, {(1,0)/e}, {(1,2)/f}, {(2,0)/g}, {(2,2)/h}}
        \node[enode] (\name) at \pos {};
    \foreach \source/\dest in {c/d,c/e,d/e,f/g,f/h,g/h}
        \path[edge] (\source) -- (\dest);
          \node[math,right] () at (2.2,1) {\Large $e$};
     \end{scope}

\end{tikzpicture}
\end{center}
\caption{Degrees 5 and 2 configuration reduction using a multiple edge}\label{f:a6red}
\end{figure}
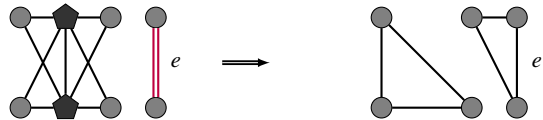
As shown in figure \ref{f:a6red}, if $c>0$ then we can add triangles until the last $c$ vertices which can then be joined to the vertices of a multiple edge $e$ in $G$ and when one edge is removed from $e$ the resulting graph will be quintic and have the triangle property.
\end{proof}

\subsection{5-Wheel subgraph}

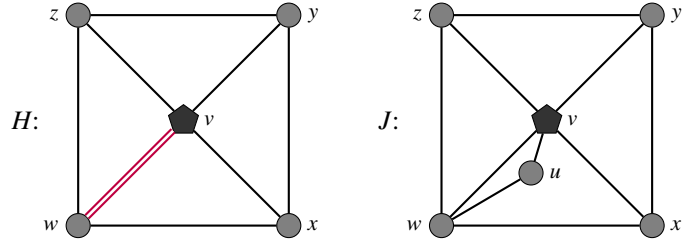
\begin{figure}[ht]
\begin{center}
\begin{tikzpicture}[scale=0.7, transform shape]
   \foreach \pos/\name in {{(0,0)/a}, {(4,0)/b}, {(4,4)/c},{(0,4)/d}}
        \node[enode] (\name) at \pos {};
   \foreach \pos/\name in {{(2,2)/v}}
        \node[pnode] (\name) at \pos {};

%        \node[nvertex] (\name) at \pos {};
    % First we draw the vertices
    % Connect vertices with edges and draw weights
    \foreach \source/\dest in {b/c,c/d,d/a,b/a,v/b,v/c,v/d}
        \path[edge] (\source) -- (\dest);
    \foreach \source/\dest in {a/v}
        \path[dedge] (\source) -- (\dest);
\node[math,right] () at (v.east) {\large $v$};
\node[math,left] () at (a.west) {\large $w$};
\node[math,right] () at (b.east) {\large $x$};
\node[math,right] () at (c.east) {\large $y$};
\node[math,left] () at (d.west) {\large $z$};
%\foreach \source/\dest in {e/f,e/e1,e/e2,e/e4,e/e5,f/f1,f/f2,f/f3,f/f4}
   %     \path[edge] (\source) -- (\dest);
\node[math] () at (-1,2) {\Large $H$:};
%\node[math] () at (5,1.5) {$\tikzLongrightarrow$};
\end{tikzpicture}
~~~~~
\begin{tikzpicture}[scale=0.7, transform shape]
   \foreach \pos/\name in {{(0,0)/a}, {(4,0)/b}, {(4,4)/c},{(0,4)/d}, {(1.7,1)/u}}
        \node[enode] (\name) at \pos {};
  \foreach \pos/\name in {{(2,2)/v}}
        \node[pnode] (\name) at \pos {};
%        \node[nvertex] (\name) at \pos {};
    % First we draw the vertices
    % Connect vertices with edges and draw weights
    \foreach \source/\dest in {a/b,b/c,c/d,d/a,v/a,v/b,v/c,v/d,v/u,a/u}
        \path[edge] (\source) -- (\dest);
\node[math] () at (-1,2) {\Large $J$:};
\node[math,right] () at (u.east) {\large $u$};
\node[math,right] () at (v.east) {\large $v$};
\node[math,left] () at (a.west) {\large $w$};
\node[math,right] () at (b.east) {\large $x$};
\node[math,right] () at (c.east) {\large $y$};
\node[math,left] () at (d.west) {\large $z$};
%\foreach \source/\dest in {e/f,e/e1,e/e2,e/e4,e/e5,f/f1,f/f2,f/f3,f/f4}
   %     \path[edge] (\source) -- (\dest);
%\node[math] () at (0,1.4) {$H$};
%\node[math] () at (5,1.5) {$\tikzLongrightarrow$};
\end{tikzpicture}
\end{center}
\caption{The base configurations for $W_5$}\label{f:w5}
\end{figure}

We will concentrate on subgraphs isomorphic to the diamond and again start by considering the case in which there is a vertex adjacent to more than two of the vertices of a diamond. Since $\omega(G)=3$ the only possibility is the wheel with five vertices, the join $W_5:=C_4 + \{v\}$. Moreover, by observation \ref{o:four}, vertex $v$ needs a fifth edge which is either a double edge to a vertex in the $C_4$ or an edge to a new vertex which is also a neighbour of a vertex in the $C_4$, as shown, without loss of generality, in figure \ref{f:w5}.

We can use two different reductions which will be able to produce smaller quintic graphs with the triangle property.
We can use X-reductions to form $H'$ using $xwvz$ or $H''$ by using $zwvx$, or, similarly, $J'$ by Z-reduction using $uvwz$ and $J''$ by using $uvwx$. Note the diagonal symmetry of these graphs as drawn in figure \ref{f:w5} means that
the results of the two different ways of reducing produce the same basic structure.

For both $H$ and $J$, vertex $w$ is adjacent to all named vertices apart from $y$ and $w$ cannot be adjacent to $y$ since that would create a $K_4$. Thus $w$'s fifth edge is either a double edge to one of its existing neighbours $x$ or $z$ (or $u$ for $J$),  or to a new neighbour of one of them to ensure it is a triangle; all $v$'s neighbours are already known. Both $v$ and $w$ are used in all of the four reductions described though.

For $H$, property X1 can hold if $xw$ is a double edge, %so $H'$ would have a loop, 
but then $H''$ is sure not to, and vice versa. Property X3 does not hold since there is a triangle involving $y$ and $v$ for both reductions. For property X2, as with the X1 case, there can only possibly be one unsafe triangle using $w$ ($v$'s triangles all have at least two vertices from $H$) and so we can choose $H'$ or $H''$ depending on whether the unsafe triangle is with $x$ or $z$.

Similarly, for $J$ we can make a similar argument based on $w$'s fifth edge. Since $w$ and $v$ are only adjacent to vertices in $J$ property Z3 cannot hold. Z1 or Z2 can only hold for $w$ with one of $x$ or $z$, so we can use whichever of $J'$ or $J''$ that does not violate the property. If $wu$ is a double edge then we can use an X-reduction and a $K_4$ is formed in both $J'$ and $J''$. 

Finally, if in $J$ there is a vertex $t$ and edges $tu$ and $tw$ that are only in one triangle then we can use a Z-reduction using $tuvw$. Again Z3 cannot hold because of $w$ and $v$. If $tu$ is a double edge, the X-reduction will work, and X2 isn't possible due to $tu$ only being in one triangle with $w$.

\subsection{No \texorpdfstring {$W_5$}{W5} subgraph}

From now on we can suppose that no diamond in $G$ has a neighbouring vertex which is adjacent to more than two of the vertices in the diamond, in addition to there being no $K_4$ subgraph. Let $t_e$ be the number of vertices
adjacent to both vertices of an edge $e$. Since $G$ has the triangle property we know that $t_e \geq 1$ for all $e\in E(G)$ and by 5-regularity we must have $t_e\leq 4$. %If $e$ is part of a triangle that uses multiple edges we can have $t_e>4$

We will deal first with the cases where we have an edge $e$ such that $t_e > 2$.

\begin{description}
\item{Case i) $t_e = 4$:}\\
Since we have four different vertices adjacent to $e$, none of these vertices can be adjacent to each other as that would form a $K_4$. Thus the only possibility is $K_{4,1,1}$ as shown on the left of figure \ref{f:a6red} and this will be another atom, $A_6$. 
None of the vertices of degree 2 in $A_6$ can be joined to each other since that would induce a $K_4$. 
%Additionally, if there is a multiple edge anywhere in the graph, we can use lemma \ref{l:2s} as in figure \ref{f:a6red}, since we have four different choices from the degree 2 vertices in $A_6$.
However, we can also choose to directly replace such an $A_6$ by $A_1$ and $A_4$; although it does not reduce the number of vertices, as mentioned in section \ref{s:intro}, all graphs with cut vertices will be proved reducible in theorem \ref{t:a4red}.

%Moreover, if any two vertices of degree 2 have a common neighbour outside of the configuration then we can reduce them as shown in figure \ref{f:a6red}. The edges to $v$ must be in a triangle, and there must be a new vertex adjacent to $u$ and $v$ as $u$'s other known neighbours are degree 5 already.

\item{Case ii) $t_e = 3$:}\\
Suppose the edge $e=u_1u_2$ as shown in figure \ref{f:te3}, and the three vertices adjacent to $e$ are $\{v_1, v_2, v_3\}$.
Since $t_e=3$, the other two neighbours of $u_1$ and $u_2$ are distinct, and suppose they are $w_1$ and $w_2$. 

\begin{figure}[ht]
\begin{center}
\begin{tikzpicture}[scale=0.6, transform shape]
  \foreach \pos/\name in {{(4,4)/v2},{(0,4)/v1}, {(2,-3)/v3}, {(-3,-1.5)/w1}, {(7,-1.5)/w2}}
        \node[enode] (\name) at \pos {};
 \foreach \pos/\name in {{(0,0)/u1}, {(4,0)/u2}}
        \node[pnode] (\name) at \pos {};
    \foreach \source/\dest in {u1/u2,u1/v1,u1/v2,u1/v3,u2/v1,u2/v2,u2/v3,u1/w1,u2/w2}
        \path[edge] (\source) -- (\dest);
\node[math,left] () at (u1.west) {\Large $u_1$};
\node[math,right] () at (u2.east) {\Large $u_2$};
\node[math,left] () at (v1.west) {\Large $v_1$};
\node[math,right] () at (v2.east) {\Large $v_2$};
\node[math,right] () at (v3.east) {\Large $v_3$};
\node[math,left] () at (w1.west) {\Large $w_1$};
\node[math,right] () at (w2.east) {\Large $w_2$};
\node[math,right] () at (2,0.3) {\Large $e$};
\end{tikzpicture}
\end{center}
\caption{General subgraph with $t_e=3$}\label{f:te3}
\end{figure}
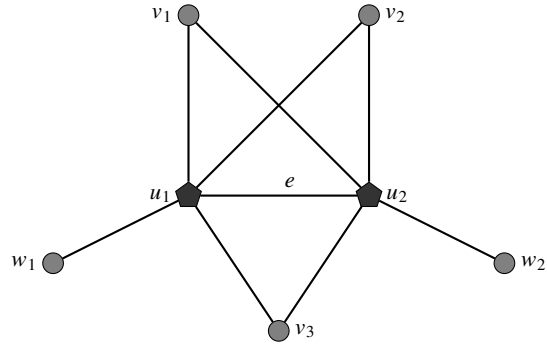

\begin{description}
\item{Subcase a) An edge incident with $e$ (or $e$ itself) is doubled:}

If we had $w_1=u_2$ (and necessarily $w_2=u_1$) we can use the X-reduction using $e$ which cannot violate any of the properties since all edges adjacent to $u_1$ and $u_2$ are known. Similarly, if $w_1=v_i$ (without loss of generality let us suppose $i=1$) then we can almost always use the X-reduction on $v_1u_2u_1v_j$ for $j\in\{2,3\}$; it will only fail if $w_2$ was either $v_1$ or was adjacent to $v_1$. However, we cannot have $w_2=v_1$ as that vertex would then have two double edges to $e$ and thus its fifth edge must be to $v_2$ or $v_3$, which would form a $K_4$. 

If $w_1=v_1$ and  $w_2v_1 \in E(G)$ the fifth edge from $v_1$ must either be to a neighbour of $w_2$ or to $w_2$. In the former case the graph resulting from the 
Z-reduction of $u_1v_1u_2w_2$ will be quintic and have the triangle property. In the latter case we know all five neighbours of each of $v_1$, $u_1$ and $u_2$ and hence $w_2$ must be adjacent to either $v_1$ with a double edge, or to two previously unnamed joined vertices, $x_1$ and $x_2$ (since if $w_2v_j\in E(G)$ then there exists a $W_5$ in $G$). We can delete $\{u_1,u_2,v_1,w_2\}$ and add edges $x_1v_3$ and $x_2v_3$ to form a triangle. $v_2$ can then be joined to two new vertices joined by a quadruple edge and the resulting graph has fewer vertices than $G$ but still is quintic and has the triangle property.

\item{Subcase b) No double edge between vertices in figure \ref{f:te3}:}

Now $w_1$ and $w_2$ are vertices not in $\{u_1,u_2,v_1,v_2,v_3\}$ and so both edges $u_jw_j$ need to be in a triangle.
Each $w_j$ can only be adjacent to at most one $v_i$ since otherwise a $W_5$ subgraph exists, contrary to our supposition.
There are thus two cases to consider; either $w_1$ and $w_2$ are joined to the same $v_i$ or different ones.

If $w_1v_3$ and $w_2v_3$ are both in $E(G)$ then we can use a Z-reduction using $v_1u_2u_1v_2$ and the resulting graph is 5-regular and has the triangle property because there are no unknown triangles from $u_1$ or $u_2$.
For the other case, if, say, $w_1v_1$ and $w_2v_2$ are both in $E(G)$ and either $v_1w_1$ or $v_2w_2$ is a multiple edge then we can use the X-reduction with it and $u_1$ and $u_2$. Note that, additionally,
 $v_1$ and $v_2$ cannot be adjacent to $v_3$ or each other as that would form a $K_4$ with $u_1$ and $u_2$, and so they were both adjacent to two previously unnamed vertices, let us call them $\{x_1, x_2, y_1, y_2\}$ where $v_kx_k$ and $v_ky_k$ were the edges in $G$.

We can now, as we did in figure \ref{f:ident23}, delete $\{u_1, u_2, v_1, v_2\}$ and add two new vertices joined by a quadruple edge joined to one of $\{v_3,w_1,w_2\}$. The other two vertices from this set can then be joined to $x_1$ and $y_1$ and $x_2$ and $y$. If $x_1y_1$ and $x_2y_2$ were edges in $G$ then all added edges are in triangles and the resulting graph is quintic.
 
However, if, say $x_1y_1$ is not an edge, then both must be joined to a neighbour of $v_1$ to have the triangle property, and $w_1$ is the only possibility for that. In this case we can delete $u_1$ and $v_1$ and add edges $\{x_1y_1, w_1u_2, w_1v_3\}$ which will give a quintic graph with the triangle property.
\end{description}
\end{description}

\subsection{For all \texorpdfstring{$e\in E(G); t_e\leq 2$}{eEte2}}

\begin{figure}[ht]
\begin{center}
\begin{tikzpicture}[scale=0.6, transform shape]
   \foreach \pos/\name in {{(0,0)/a}, {(4,0)/b}, {(2,-2)/c},{(2,2)/d}}
        \node[enode] (\name) at \pos {};
   \foreach \pos/\name in {{(2,-2)/c},{(2,2)/d}}
        \node[pnode] (\name) at \pos {};
    \foreach \source/\dest in {a/c, a/d,b/c,b/d}
        \path[edge] (\source) -- (\dest);
    \foreach \source/\dest in {c/d}
        \path[tredge] (\source) -- (\dest);
\node[math] () at (-1,1.5) {\LARGE $A_7:$};
\node[math,left] () at (a.west) {\Large $u_1$};
\node[math,right] () at (b.east) {\Large $u_2$};

 \begin{scope}[shift={(5,-2)}]
  \foreach \pos/\name in {{(4,0)/d}, {(5,4)/e}, {(6,0)/f}}
        \node[enode] (\name) at \pos {};
  \foreach \pos/\name in {{(2,0)/a}, {(1,4)/b}, {(3,4)/c}}
        \node[pnode] (\name) at \pos {};
    \foreach \source/\dest in {a/d,a/c,c/d,c/e,d/e,d/f,e/f}
        \path[edge] (\source) -- (\dest);
    \foreach \source/\dest in {a/b}
        \path[tredge] (\source) -- (\dest);
    \foreach \source/\dest in {c/b}
        \path[dedge] (\source) -- (\dest);
\node[math,below] () at (d.south) {\Large $v_1$};
\node[math,right] () at (e.east) {\Large $v_2$};
\node[math,right] () at (f.east) {\Large $v_3$};
\end{scope}
\begin{scope}[shift={(12,-2)}]
  \foreach \pos/\name in {{(3,4)/c},{(3,2)/d}, {(3,0)/e}}
        \node[enode] (\name) at \pos {};
  \foreach \pos/\name in {{(1,0)/a}, {(1,4)/b}}
        \node[pnode] (\name) at \pos {};
    \foreach \source/\dest in {a/d,a/e,b/d,b/c,c/d,d/e}
        \path[edge] (\source) -- (\dest);
    \foreach \source/\dest in {a/b}
        \path[tredge] (\source) -- (\dest);
\node[math,right] () at (c.east) {\Large $w_1$};
\node[math,right] () at (d.east) {\Large $w_2$};
\node[math,right] () at (e.east) {\Large $w_3$};
\end{scope}
\end{tikzpicture}
\end{center}
\caption{Triple edges in subgraphs with \texorpdfstring{$t_e = 2$}{te2}}\label{f:trip}
\end{figure}
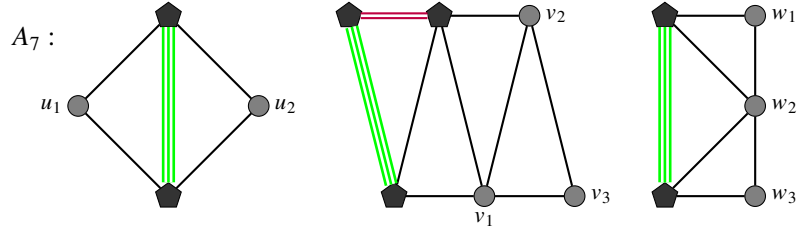

Firstly, we can deal with multiple edges which have are more than doubled; if there is a quadruple edge then it can only be attached to a cut vertex as in figure \ref{f:cutv}. For a triple edge there will also be a cut-vertex unless one of the subgraphs shown in figure \ref{f:trip}
exists. We will now deal with these in sequence, letting $H$ be each of the three different configurations; if $u_1$ and $u_2$ have a common neighbour then we can use the X-reduction, otherwise this configuration is the atom $A_7$. By lemma \ref{l:2s} we can reduce $A_7$ if there are two multiple edges anywhere in $G$ outside of $H$.

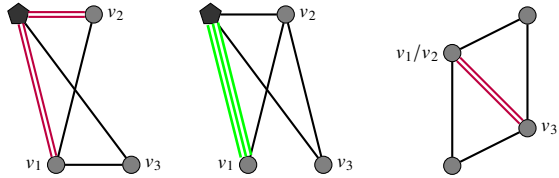
\begin{figure}[ht]
\begin{center}
\begin{tikzpicture}[scale=0.5, transform shape]
  \foreach \pos/\name in {{(5,4)/c},{(4,0)/e}, {(6,0)/f}}
        \node[enode] (\name) at \pos {};
  \foreach \pos/\name in {{(3,4)/d}}
        \node[pnode] (\name) at \pos {};
    \foreach \source/\dest in {f/d,f/e,e/c}
        \path[edge] (\source) -- (\dest);
    \foreach \source/\dest in {c/d,d/e}
        \path[dedge] (\source) -- (\dest);
\node[math,right] () at (c.east) {\Large $v_2$};
\node[math,left] () at (e.west) {\Large $v_1$};
\node[math,right] () at (f.east) {\Large $v_3$};
\end{tikzpicture}
~~~
\begin{tikzpicture}[scale=0.5, transform shape]
  \foreach \pos/\name in {{(4,0)/c},{(5,4)/e}, {(6,0)/f}}
        \node[enode] (\name) at \pos {};
  \foreach \pos/\name in {{(3,4)/d}}
        \node[pnode] (\name) at \pos {};
   \foreach \source/\dest in {f/d,f/e,e/c,d/e}
        \path[edge] (\source) -- (\dest);
    \foreach \source/\dest in {c/d}
        \path[tredge] (\source) -- (\dest);
\node[math,left] () at (c.west) {\Large $v_1$};
\node[math,right] () at (e.east) {\Large $v_2$};
\node[math,right] () at (f.east) {\Large $v_3$};
\end{tikzpicture}
~~~
\begin{tikzpicture}[scale=0.5, transform shape]
  \foreach \pos/\name in {{(4,2)/a}, {(6,3)/b}, {(4,-1)/e}, {(6,0)/f}}
        \node[enode] (\name) at \pos {};
    \foreach \source/\dest in {a/b,b/f,e/f,e/a}
        \path[edge] (\source) -- (\dest);
    \foreach \source/\dest in {f/a}
        \path[dedge] (\source) -- (\dest);
        \node[math,left] () at (a.west) {\Large $v_1/v_2$};
\node[math,right] () at (f.east) {\Large $v_3$};
\end{tikzpicture}

%\begin{tikzpicture}[scale=0.5, transform shape]
%  \foreach \pos/\name in {{(2,0)/a}, {(1,4)/b}, {(3,4)/c},{(4,0)/d}, {(5,4)/e}, {(6,0)/f}}
%        \node[pnode] (\name) at \pos {};
%    \foreach \source/\dest in {a/d,a/c,c/d,c/e,d/e,d/f}
%        \path[edge] (\source) -- (\dest);
%    \foreach \source/\dest in {a/b,e/f}
%        \path[tredge] (\source) -- (\dest);
%    \foreach \source/\dest in {c/b,f/d}
%        \path[dedge] (\source) -- (\dest);
%\end{tikzpicture}
\end{center}
\caption{The three reductions necessary for the second configuration}\label{f:v1v2}
\end{figure}

In the second configuration we can use the two possible reductions shown on the left of figure \ref{f:v1v2} which are quintic since the degrees of $v_1$, $v_2$ and $v_3$ remain 4, 3 and 2 (respectively) in the configuration. The left reduction will give a graph with the triangle property unless $v_2v_3$ was part of a triangle $T$ in $G$ and was not a multiple edge and the middle reduction will require that $v_1v_3$ is part of a similar triangle.%, and so for neither of the reductions to give graphs with the triangle property, 

Thus, $G$ is only not reducible to a graph with the triangle property by one of these two operations if both $v_2v_3$ and $v_1v_3$ are non multiple edges in triangles with new different vertices. However, in that case we can contract $v_1$, $v_2$ and the three pentagonal vertices into one vertex as shown in the right of figure \ref{f:v1v2} and the resulting graph will be quintic and have the triangle property.

% If both $v_2v_3$ and $v_1v_3$ are multiple edges, then the former must be a triple edge and we get another irreducible 5-regular graph with the triangle property and six vertices shown on the right of figure \ref{f:v1v2}.

We can argue similarly for the right hand graph in figure \ref{f:trip}. Vertex $w_2$'s fifth edge can either be (without loss of generality) to $w_1$, or it can be 
part of a triangle with $w_1$ and/or $w_3$. Note that $w_1w_3$ cannot be an edge since we are assuming $G$ does not contain $W_5$. 
If $w_1w_2$ is a double edge then we can contract all five vertices in the configuration into one and the reduced graph will still have the triangle property. If the fifth edge from $w_2$ is in a triangle with both $w_1$ and $w_3$ then we can use the X-reduction centred on the triple edge.

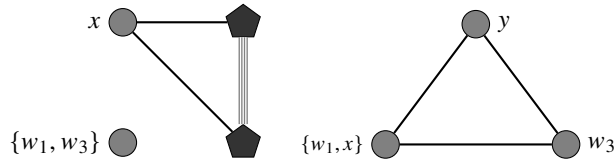
\begin{figure}[ht]
\begin{center}
\begin{tikzpicture}[scale=0.8, transform shape]
  \foreach \pos/\name in {{(0,2)/uv},{(0,4)/w}}
        \node[enode] (\name) at \pos {};
  \foreach \pos/\name in {{(2,4)/d}, {(2,2)/e}}
        \node[pnode] (\name) at \pos {};
    \foreach \source/\dest in {w/d,w/e}
        \path[edge] (\source) -- (\dest);
    \foreach \source/\dest in {e/d}
        \path[quple] (\source) -- (\dest);
\node[math,left] () at (uv.west) {\large $\{w_1,w_3\}$};
\node[math,left] () at (w.west) {\large $x$};
\end{tikzpicture}
~~~
\begin{tikzpicture}[scale=0.8, transform shape]
  \foreach \pos/\name in {{(0,0)/c},{(3,0)/d}, {(1.5,2)/e}}
        \node[enode] (\name) at \pos {};
    \foreach \source/\dest in {c/d,d/e,e/c}
        \path[edge] (\source) -- (\dest);
\node[math,left] () at (c.west) {$\{w_1,x\}$};
\node[math,right] () at (d.east) {\large $w_3$};
\node[math,right] () at (e.east) {\large $y$};
\end{tikzpicture}
\end{center}
\caption{Two possible reductions for the right graph in figure \ref{f:trip}}\label{f:w1w2}
\end{figure}

Lastly, we can suppose there is a vertex $x$ adjacent to both $w_1$ and $w_2$, but not $w_3$. Again there are two possible reductions we can use, as shown in figure \ref{f:w1w2}; we can either identify the five vertices as we did when $w_1w_2$ was a double edge and add a quadruple edge pendant to $x$, or (if $w_1x$ is part of a triangle with both $w_2$ and a new vertex $y$) contract $w_1x$ and rejoin as shown.

\section{Double Edges are the only multiple edges}
\begin{figure}[ht]
\begin{center}
\begin{tikzpicture}[scale=0.7, transform shape]
   \foreach \pos/\name in {{(0,1)/a}, {(1,0)/b}, {(1,2)/c}}
        \node[enode] (\name) at \pos {};
\foreach \pos/\name in {{(-1,1)/a1}, {(-2,0)/b1}, {(-2,2)/c1}, {(-3,1)/d1}}
        \node[enode] (\name) at \pos {};
     \foreach \pos/\name in {{(-0.5,2)/a2}, {(-1.5,3)/b2}, {(0.5,3)/c2}}
        \node[enode] (\name) at \pos {};
      \foreach \pos/\name in {{(-0.5,0)/a3}, {(-1.5,-1)/b3}, {(0.5,-1)/c3}}
        \node[enode] (\name) at \pos {};
    \foreach \source/\dest in {a/b,b/c,c/a,a1/b1,a1/c1,d1/b1,d1/c1,b1/c1,a2/b2,a2/c2,b2/c2,a3/b3,a3/c3,b3/c3}
        \path[edge] (\source) -- (\dest);
\draw (-0.5,1) ellipse (0.9cm and 1.3cm);

\node[math] () at (3,1) {$\tikzLongrightarrow$};

   \foreach \pos/\name in {{(8.5,0)/a}, {(8.5,2)/b}, {(9.5,0)/c}, {(9.5,2)/d},{(6.5,1)/e},{(7.5,1)/f}}
        \node[enode] (\name) at \pos {};
 \foreach \pos/\name in {{(5.5,0)/b1}, {(5.5,2)/c1}, {(4.5,1)/d1}}
        \node[enode] (\name) at \pos {};
         \foreach \pos/\name in {{(6.5,3)/b2}, {(7.5,3)/c2}}
        \node[enode] (\name) at \pos {};
      \foreach \pos/\name in {{(6.5,-1)/b3}, {(7.5,-1)/c3}}
        \node[enode] (\name) at \pos {};
 \foreach \source/\dest in {a/c,b/c,a/d,b/d,d1/b1,d1/c1,e/f,e/b2,e/b3,f/b2,f/b3,e/c2,e/c3,f/c2,f/c3}
        \path[edge] (\source) -- (\dest);
\foreach \source/\dest in {a/b}
        \path[tredge] (\source) -- (\dest);
 \foreach \source/\dest in {b1/c1}
        \path[dedge] (\source) -- (\dest);
\end{tikzpicture}
\end{center}
\caption{Pendant triangle reduction options}\label{f:trigo}
\end{figure}
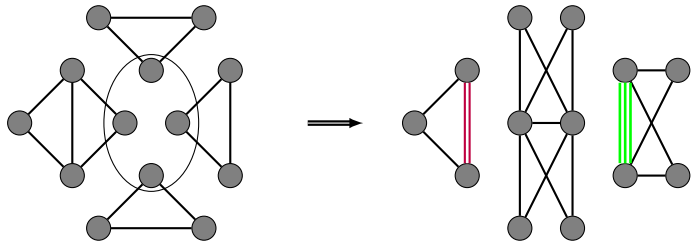

Since we have either reduced or shown atomic %or irreducible 
all graphs in this section with triple or quadruple edges, we can now assume any multiple edges are double edges. %Recall that triangles which are edge disjoint from a configuration but have one vertex in it are called pendant triangles. 
The edge between the other two vertices in a pendant triangle may or may not be needed in the reduced graph to ensure the triangle property holds for other edges. To address this we can again use the reductions in figure \ref{f:ptri}; note that when introducing triple edges we will increase the number of vertices, as on the right in figure \ref{f:trigo}, but this will still result in a reduction so long as the configuration in the oval which was deleted had sufficient vertices.

\subsection{Aloof triangles}
If any aloof triangle in $G$ contains more than one double edge then one vertex $v$ in the triangle must be adjacent to both double edges and hence, by observation \ref{o:four} $v$'s fifth edge must be in a triangle with one of the other neighbours of $v$, forming a diamond, a contradiction of aloofness. Similarly, if all edges of a triangle were double edges then we must have a $K_4$ and a cut-vertex as shown in figure \ref{f:cutv}.

\begin{figure}[ht]
\begin{center}
\begin{tikzpicture}[scale=0.5, transform shape]
    \begin{scope}[shift={(-1,0)}]
     \foreach \pos/\name in {{(0,2.7)/c},{(-1,1)/e},{(1,1)/f}}
        \node[enode] (\name) at \pos {};
     \foreach \source/\dest in {f/c,f/e}
        \path[edge] (\source) -- (\dest);
    \foreach \source/\dest in {c/e}
        \path[dedge] (\source) -- (\dest);
 \node[math] () at (2,2) {\Large :};
 \node[math,left] () at (-1,1.5) {{\Large $A_3$:}};
 \end{scope}
 
  \foreach \pos/\name in { {(4.5,4)/c1},{(3,1)/e1}, {(5.5,4)/c2},{(3.5,0)/e2},{(2.5,0)/e3}, {(6,1)/f}}
        \node[enode] (\name) at \pos {};
  \foreach \pos/\name in { {(5,2.7)/c},{(4,1)/e}}
        \node[pnode] (\name) at \pos {};
    \foreach \source/\dest in {f/c,f/e,c/c1,c/c2,e/e1,e/e2,e1/e2,c1/c2,e1/e3,e2/e3}
        \path[edge] (\source) -- (\dest);
    \foreach \source/\dest in {c/e}
        \path[dedge] (\source) -- (\dest);

\node[math] () at (7.5,2) {$\tikzLongrightarrow$};

 \foreach \pos/\name in {{(10.5,4)/c1},{(9,1)/e1}, {(11.5,4)/c2},{(9.5,0)/e2},{(8.5,0)/e3}, {(12,1)/f}}
        \node[enode] (\name) at \pos {};
    \foreach \source/\dest in {c1/f,c2/f,c1/c2,e1/e3,e2/e3}
        \path[edge] (\source) -- (\dest);
   \foreach \source/\dest in {e1/e2}
        \path[dedge] (\source) -- (\dest);
        
%        \begin{scope}[shift={(0,-6)}]
%
%  \foreach \pos/\name in { {(4.5,4)/c1},{(3,1)/e1}, {(5.5,4)/c2},{(3.5,0)/e2},{(3.25,2.5)/e3}, {(6,1)/f}}
%        \node[enode] (\name) at \pos {};
%  \foreach \pos/\name in { {(5,2.7)/c},{(4,1)/e}}
%        \node[pnode] (\name) at \pos {};
%    \foreach \source/\dest in {f/c,f/e,c/c1,c/c2,e/e1,e/e2,e1/e2,c1/c2,e1/e3,c1/e3}
%        \path[edge] (\source) -- (\dest);
%    \foreach \source/\dest in {c/e}
%        \path[dedge] (\source) -- (\dest);
%\node[math,left] () at (e3.west) {{\Large $y$}};
%\node[math] () at (7.5,2) {$\tikzLongrightarrow$};
%
% \foreach \pos/\name in {{(10.5,4)/c1},{(9,1)/e1}, {(11.5,4)/c2},{(9.5,0)/e2},{(9.25,2.5)/e3}, {(12,1)/f}}
%        \node[enode] (\name) at \pos {};
%    \foreach \source/\dest in {e2/f,c2/f,e2/c2,e1/e3,c1/e3,c1/e1}
%        \path[edge] (\source) -- (\dest);
%        \node[math,left] () at (e3.west) {{\Large $y$}};
%        \end{scope}
\end{tikzpicture}
\end{center}
\caption{Lone triangle with double edge}\label{f:dbintri}
\end{figure}
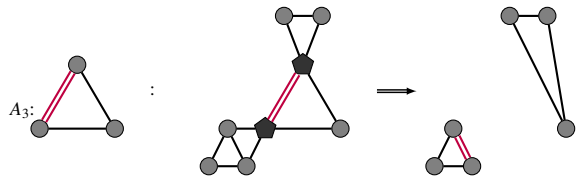

If we have an aloof triangle $T$ which contains exactly one double edge this is atom $A_3$ which is shown in figure \ref{f:dbintri}. If $T$ is not part of a diamond then the neighbours of each vertex in $T$ must be disjoint, as in the top of figure \ref{f:dbintri}, and if not both pendant triangles from the degree 3 vertices of $T$ are aloof then we can delete the vertices joined by the double edge and form the reduced graph as shown. 
Only if both triangles are aloof is $T$ not reducible in this way, similar to atom $A_2$.
%% but it is reducible if there are two multiple edges elsewhere in the graph, by lemma \ref{l:2s}. 
%Additionally, if any two of the five vertices of degree 2 have a common neighbour outside of the configuration, as shown with vertex $y$ in the bottom of figure \ref{f:dbintri} it is also reducible by joining two neighbours of $y$ and then making the other three vertices all adjacent.
%%and we will deal with this in section \ref{s:atom} too.
%%However, in this case, if there are two multiple edges anywhere else in $G$ then we can delete both of the pentagonal vertices of $T$, similarly to what was done for $A_8$ in figure \ref{f:a6red}.

%Henceforth, we can assume that any double edge must be in a diamond.

\subsection{Three Double Edges in a diamond}
\begin{figure}[ht]
\begin{center}
\begin{tikzpicture}[scale=0.7, transform shape]
  \foreach \pos/\name in {{(3,4)/c},{(5,4)/e}}
        \node[enode] (\name) at \pos {};
  \foreach \pos/\name in {{(4,2.5)/f}, {(4,5.5)/g}}
        \node[pnode] (\name) at \pos {};
   \foreach \source/\dest in {f/c,g/e}
        \path[edge] (\source) -- (\dest);
    \foreach \source/\dest in {c/g,e/f,f/g}
        \path[dedge] (\source) -- (\dest);
\node[math,left] () at (2,4) {{\Large $A_8$:}};
\end{tikzpicture}
~~~~,~~~~~~~~
\begin{tikzpicture}[scale=0.7, transform shape]
  \foreach \pos/\name in {{(3,4)/c},{(5,4)/e}, {(4,2.5)/f}}
        \node[enode] (\name) at \pos {};
  \foreach \pos/\name in {{(4,5.5)/g}}
        \node[pnode] (\name) at \pos {};
   \foreach \source/\dest in {g/c,g/e,f/e}
        \path[dedge] (\source) -- (\dest);
    \foreach \source/\dest in {f/g,f/c}
        \path[edge] (\source) -- (\dest);
\node[math,right] () at (e.east) {{\Large $v$}};
\node[math,right] () at (f.east) {{\Large $w$}};
\end{tikzpicture}
\end{center}
\caption{Three double edges in a diamond}\label{f:threedoub}
\end{figure}
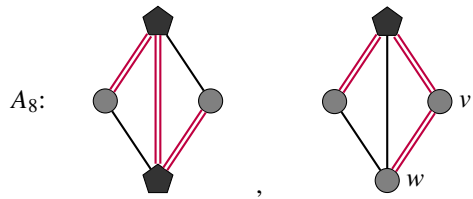

It is possible to have three double edges in a diamond in two ways as shown in figure \ref{f:threedoub}.
In the left case ($A_8$) we consider the neighbours of the grey vertices to find a reduction; if they share a neighbour then we can use the X-reduction to delete the degree 5 vertices and add a triple edge between the grey vertices so the triangle property still holds. If they do not share a neighbour then there is a pendant triangle from each grey vertex as shown in figure \ref{f:trigo}. After deleting all the vertices of $A_8$ we have two choices: as in lemma \ref{l:a5red} we can double the edges of the unsafe triangles and add two new vertices for the aloof triangles. In each case the resulting graph will be 5-regular and have the triangle property, and will have fewer vertices, so configuration $A_8$ is reducible. 
%unless all its neighbours are in pendant aloof triangles, and, moreover, by lemma \ref{l:2s} the whole graph cannot have any other multiple edges.

For the right configuration in figure \ref{f:threedoub} vertex $v$ needs a fifth edge, and in order to be in a triangle without creating a triple edge 
it must be in a new triangle adjacent to $w$. However, this gives a configuration with three vertices of degree 5 and one each of degree 2 and 3; this can be reduced to a single vertex since there cannot be any other edges between the vertices without there being a previously dealt with subgraph.

\subsection{Two Double Edges in a diamond}

\begin{figure}[ht]
\begin{center}
\begin{tikzpicture}[scale=0.55, transform shape]
  \foreach \pos/\name in {{(3,4)/c},{(4,3)/f},{(5,4)/e}}
        \node[enode] (\name) at \pos {};
  \foreach \pos/\name in {{(4,5)/g}}
        \node[pnode] (\name) at \pos {};
   \foreach \source/\dest in {f/c,g/e,f/e}
        \path[edge] (\source) -- (\dest);
    \foreach \source/\dest in {c/g,f/g}
        \path[dedge] (\source) -- (\dest);
\node[math,right] () at (f.east) {{\LARGE $v$}};

\node[math,rotate=22.5] () at (6,5) {\Large :};

  \foreach \pos/\name in {{(7,6)/c},{(9,4)/d},{(9,6)/e}}
        \node[enode] (\name) at \pos {};
  \foreach \pos/\name in {{(8,5)/f},{(8,7)/g}}
        \node[pnode] (\name) at \pos {};
   \foreach \source/\dest in {f/c,g/e,f/e,d/e,d/f}
        \path[edge] (\source) -- (\dest);
    \foreach \source/\dest in {c/g,f/g}
        \path[dedge] (\source) -- (\dest);

\node[math,rotate=-22.5] () at (6,3) {\Large :};
\node[math,rotate=0] () at (6,4) {\Large or};

  \foreach \pos/\name in {{(7,2)/c},{(7,0)/d},{(9,2)/e}}
        \node[enode] (\name) at \pos {};
  \foreach \pos/\name in {{(8,1)/f},{(8,3)/g}}
        \node[pnode] (\name) at \pos {};
   \foreach \source/\dest in {f/c,g/e,f/e,d/c,d/f}
        \path[edge] (\source) -- (\dest);
    \foreach \source/\dest in {c/g,f/g}
        \path[dedge] (\source) -- (\dest);

\node[math] () at (11,5) {$\tikzLongrightarrow$};

  \foreach \pos/\name in {{(12,6)/c},{(14,4)/d},{(14,6)/e}}
        \node[enode] (\name) at \pos {};
     \foreach \source/\dest in {d/c,d/e}
        \path[edge] (\source) -- (\dest);
    \foreach \source/\dest in {c/e}
        \path[dedge] (\source) -- (\dest);

\node[math] () at (10,1.5) {\Large :};

  \foreach \pos/\name in {{(11,1)/b},{(12,0)/d},{(14,2)/e}}
        \node[enode] (\name) at \pos {};
  \foreach \pos/\name in {{(12,2)/c},{(13,1)/f},{(13,3)/g}}
        \node[pnode] (\name) at \pos {};
   \foreach \source/\dest in {f/c,g/e,f/e,d/c,d/f,b/c,b/d}
        \path[edge] (\source) -- (\dest);
    \foreach \source/\dest in {c/g,f/g}
        \path[dedge] (\source) -- (\dest);

\node[math] () at (15,1.5) {$\tikzLongrightarrow$};
\begin{scope}[shift={(1,0)}]
  \foreach \pos/\name in {{(15,1)/b},{(18,2)/e}}
        \node[enode] (\name) at \pos {};
  \foreach \pos/\name in {{(16,2)/c},{(17,1)/f}}
        \node[pnode] (\name) at \pos {};
   \foreach \source/\dest in {b/c,b/f}
        \path[edge] (\source) -- (\dest);
    \foreach \source/\dest in {c/f}
        \path[quple] (\source) -- (\dest);
\end{scope}
\end{tikzpicture}
\end{center}
\caption{Two incident double edges in a diamond, reduced}\label{f:twodoubi}
\end{figure}
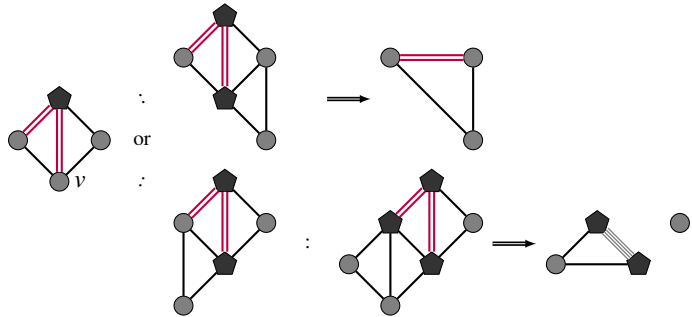

If, in a diamond, there are two double edges, they can either have a vertex in common (in two ways) or not. Let us deal with one of the former cases as shown in figure \ref{f:twodoubi}. We shall be using observation \ref{o:four} here, so that there must be a new vertex in a triangle with $v$ and one of its neighbours, and this gives us two possible configurations. Note that the lower of these two configurations must also result from when a diamond contains two double edges which meet at a vertex and still need one more edge to be added.

However, in both of these cases we can deduce the structure at the vertices of degree 4 by assuming none of the already dealt with structures exists and then reduce as shown. The resulting graphs will necessarily be quintic, have the triangle property and have fewer vertices as required.

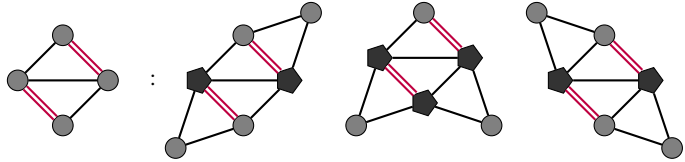
\begin{figure}[ht]
\begin{center}
\begin{tikzpicture}[scale=0.6,rotate=-90, transform shape]
  \foreach \pos/\name in {{(9,-3)/g},{(8,-4)/c},{(9,-5)/f},{(10,-4)/e}}
        \node[enode] (\name) at \pos {};
   \foreach \source/\dest in {f/c,g/e,f/g}
        \path[edge] (\source) -- (\dest);
    \foreach \source/\dest in {c/g,f/e}
        \path[dedge] (\source) -- (\dest);
%\node[math,right] () at (f.east) {{$v$}};
\node[math,rotate=90] () at (9,-2) {\Large :};

%\node[math,rotate=45] () at (6,6) {$\tikzLongrightarrow$};

  \foreach \pos/\name in {{(8,8)/c},{(10,8)/e},{(10.5,9.5)/b},{(7.5,6.5)/a}}
        \node[enode] (\name) at \pos {};
\foreach \pos/\name in {{(9,9)/g},{(9,7)/f}}
        \node[pnode] (\name) at \pos {};
   \foreach \source/\dest in {f/c,g/e,f/g,b/e,b/g,a/c,a/f}
        \path[edge] (\source) -- (\dest);
    \foreach \source/\dest in {c/g,f/e}
        \path[dedge] (\source) -- (\dest);

%\node[math,rotate=0] () at (6.25,4) {$\tikzLongrightarrow$};

  \foreach \pos/\name in {{(7.5,4)/c},{(10,2.5)/a},{(10,5.5)/b}}
        \node[enode] (\name) at \pos {};
\foreach \pos/\name in {{(8.5,5)/g},{(9.5,4)/e},{(8.5,3)/f}}
        \node[pnode] (\name) at \pos {};
   \foreach \source/\dest in {f/c,g/e,f/g,b/e,b/g,a/e,a/f}
        \path[edge] (\source) -- (\dest);
    \foreach \source/\dest in {c/g,f/e}
        \path[dedge] (\source) -- (\dest);

%\node[math,rotate=-45] () at (6,2) {$\tikzLongrightarrow$};

  \foreach \pos/\name in {{(8,0)/c},{(10,0)/e},{(10.5,-1.5)/b},{(7.5,1.5)/a}}
        \node[enode] (\name) at \pos {};
\foreach \pos/\name in {{(9,1)/g},{(9,-1)/f}}
        \node[pnode] (\name) at \pos {};
   \foreach \source/\dest in {f/c,g/e,f/g,b/e,b/f,a/c,a/g}
        \path[edge] (\source) -- (\dest);
    \foreach \source/\dest in {c/g,f/e}
        \path[dedge] (\source) -- (\dest);
\end{tikzpicture}
\end{center}
\caption{Two parallel double edges, and their structures}\label{f:pardoubi}
\end{figure}

Secondly, the two double edges in the diamond can have no vertices in common as in figure \ref{f:pardoubi}. The two vertices of degree 4 in the left configuration need a triangle with a new vertex added to them, and this can be from the three different structures shown to its right. Either the triangles are added to neighbouring or opposite edges of the diamond, but the second and fourth structures can be reduced using the X-reduction, leaving a double edge in the centre of a diamond, and the third can be treated analagously to the bottom case in figure \ref{f:twodoubi}.

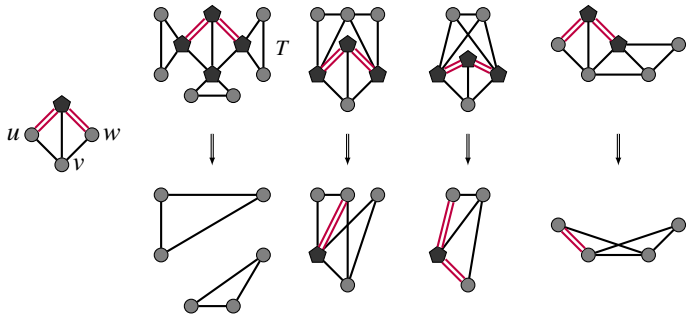
\begin{figure}[ht]
\begin{center}
\begin{tikzpicture}[scale=0.4,rotate=0, transform shape]
  \foreach \pos/\name in {{(3,4)/c},{(4,3)/f},{(5,4)/e}}
        \node[enode] (\name) at \pos {};
 \foreach \pos/\name in {{(4,5)/g}}
        \node[pnode] (\name) at \pos {};
   \foreach \source/\dest in {f/c,f/e,f/g}
        \path[edge] (\source) -- (\dest);
    \foreach \source/\dest in {c/g,g/e}
        \path[dedge] (\source) -- (\dest);
\node[math,left] () at (c.west) {{\Huge $u$}};
\node[math,right] () at (f.east) {{\Huge $v$}};
\node[math,right] () at (e.east) {{\Huge $w$}};

  \foreach \pos/\name in {{(7.3,6)/c1},{(7.3,8)/c2},{(9.7,5.3)/f1},{(8.3,5.3)/f2},{(10.7,6)/e1},{(10.7,8)/e2}}
        \node[enode] (\name) at \pos {};
 \foreach \pos/\name in {{(9,8)/g},{(8,7)/c},{(9,6)/f},{(10,7)/e}}
        \node[pnode] (\name) at \pos {};
   \foreach \source/\dest in {f/c,f/e,f/g,c/c1,c/c2,c1/c2,e/e1,e/e2,e1/e2,f/f1,f/f2,f1/f2}
        \path[edge] (\source) -- (\dest);
    \foreach \source/\dest in {c/g,g/e}
        \path[dedge] (\source) -- (\dest);
        \node[math,left] () at (11.7,6.9) {{\huge $T$}};

  \foreach \pos/\name in {{(13.5,5)/f},{(12.5,8)/h},{(13.5,8)/i},{(14.5,8)/j}}
        \node[enode] (\name) at \pos {};
 \foreach \pos/\name in {{(13.5,7)/g},{(12.5,6)/c},{(14.5,6)/e}}
        \node[pnode] (\name) at \pos {};
   \foreach \source/\dest in {f/c,f/e,f/g,c/h,c/i,h/i,e/i,e/j,i/j}
        \path[edge] (\source) -- (\dest);
    \foreach \source/\dest in {c/g,g/e}
        \path[dedge] (\source) -- (\dest);

  \foreach \pos/\name in {{(17.5,5)/f},{(17,8)/h},{(18,8)/i}}
        \node[enode] (\name) at \pos {};
 \foreach \pos/\name in {{(17.5,6.5)/g},{(16.5,6)/c},{(18.5,6)/e}}
        \node[pnode] (\name) at \pos {};
   \foreach \source/\dest in {f/c,f/e,f/g,c/h,c/i,h/i,e/i,e/h}
        \path[edge] (\source) -- (\dest);
    \foreach \source/\dest in {c/g,g/e}
        \path[dedge] (\source) -- (\dest);

  \foreach \pos/\name in {{(20.5,7)/c},{(21.5,6)/f},{(23.5,6)/f1},{(24.5,7)/e1}}
        \node[enode] (\name) at \pos {};
\foreach \pos/\name in {{(21.5,8)/g},{(22.5,7)/e}}
        \node[pnode] (\name) at \pos {};
   \foreach \source/\dest in {f/c,f/e,f/g,f/f1,e/e1,e1/f1,e/f1}
        \path[edge] (\source) -- (\dest);
    \foreach \source/\dest in {c/g,g/e}
        \path[dedge] (\source) -- (\dest);

\node[math,rotate=-90] () at (9,3.5) {$\tikzLongrightarrow$};
\node[math,rotate=-90] () at (13.5,3.5) {$\tikzLongrightarrow$};
\node[math,rotate=-90] () at (17.5,3.5) {$\tikzLongrightarrow$};
\node[math,rotate=-90] () at (22.5,3.5) {$\tikzLongrightarrow$};
% reduce to:

  \foreach \pos/\name in {{(7.3,0)/c1},{(7.3,2)/c2},{(9.7,-1.7)/f1},{(8.3,-1.7)/f2},{(10.7,0)/e1},{(10.7,2)/e2}}
        \node[enode] (\name) at \pos {};
   \foreach \source/\dest in {c1/c2,f1/f2,c1/e2,c2/e2,f1/e1,f2/e1}
        \path[edge] (\source) -- (\dest);

  \foreach \pos/\name in {{(13.5,-1)/f},{(12.5,2)/h},{(13.5,2)/i},{(14.5,2)/j}}
        \node[enode] (\name) at \pos {};
 \foreach \pos/\name in {{(12.5,0)/c}}
        \node[pnode] (\name) at \pos {};
   \foreach \source/\dest in {c/h,h/i,i/f,f/j,j/c,c/f}
        \path[edge] (\source) -- (\dest);
    \foreach \source/\dest in {c/i}
        \path[dedge] (\source) -- (\dest);

  \foreach \pos/\name in {{(17.5,-1)/f},{(17,2)/h},{(18,2)/i}}
        \node[enode] (\name) at \pos {};
 \foreach \pos/\name in {{(16.5,0)/c}}
        \node[pnode] (\name) at \pos {};
   \foreach \source/\dest in {f/i,c/i,h/i}
        \path[edge] (\source) -- (\dest);
    \foreach \source/\dest in {f/c,c/h}
        \path[dedge] (\source) -- (\dest);

  \foreach \pos/\name in {{(20.5,1)/c},{(21.5,0)/f},{(23.5,0)/f1},{(24.5,1)/e1}}
        \node[enode] (\name) at \pos {};
   \foreach \source/\dest in {f/f1,f/e1,c/f1,e1/f1}
        \path[edge] (\source) -- (\dest);
    \foreach \source/\dest in {f/c}
        \path[dedge] (\source) -- (\dest);

\end{tikzpicture}
\end{center}
\caption{Diamond with double edges and its other reductions}\label{f:atomdoub}
\end{figure}

Finally, there is one more configuration of a diamond with two double edges that is shown in figure \ref{f:atomdoub}.
If any of the vertices in the configuration are incident with disjoint pendant simple aloof triangles then it is possible to reduce using lemma \ref{l:2s} as shown for triangle $T$
in the second column of figure \ref{f:atomdoub}.
If none of the pendant triangles are aloof then we can delete the whole configuration and double the remaining edges from the triangles that were removed which will still have the triangle property and be quintic.
%, and if these triangles are not simple then they are $A_2$ or $A_3$ which allow the graph to be reduced.

If we do not have three pendant triangles then some pair of named vertices must have a new vertex in common. By symmetry we can consider $w$ which needs to have two more edges added to it, and either it can have one or two new common neighbours with $u$ or one new neighbour with $v$. It cannot have two as no edge can be in more than two triangles, and if only one then the new vertex must be in two triangles with $w$.
These three structures are shown on the right of figure \ref{f:atomdoub}, along with their reductions.

\subsection{One Double Edge in a diamond}

\begin{figure}[ht]
\begin{center}
\begin{tikzpicture}[scale=0.35, transform shape]
\begin{scope}[shift={(0,-2)}]
  \foreach \pos/\name in {{(3,6)/b},{(3,4)/c},{(4,3)/f},{(5,4)/e}}
        \node[enode] (\name) at \pos {};
  \foreach \pos/\name in {{(4,5)/g}}
        \node[pnode] (\name) at \pos {};
   \foreach \source/\dest in {f/c,g/f,f/e,c/b,c/g,g/b}
        \path[edge] (\source) -- (\dest);
    \foreach \source/\dest in {e/g}
        \path[dedge] (\source) -- (\dest);
\node[math,right] () at (e.east) {{\LARGE $v$}};
\node[math,right] () at (f.east) {{\LARGE $w$}};
\node[math,right] () at (g.east) {{\LARGE $u$}};
\node[math,left] () at (c.west) {{\LARGE $x$}};
\node[math,left] () at (b.west) {{\LARGE $y$}};
\end{scope}

\node[math,rotate=0] () at (6,2) {\Large :};

\begin{scope}[shift={(4,2)}]
\foreach \pos/\name in {{(3,6)/b},{(3,4)/c},{(4,3)/f},{(6,3)/e1},{(7,4)/e2}}
        \node[enode] (\name) at \pos {};
  \foreach \pos/\name in {{(4,5)/g},{(5,4)/e}}
        \node[pnode] (\name) at \pos {};
   \foreach \source/\dest in {f/c,g/f,f/e,c/b,c/g,g/b,f/e1,e1/e2,e/e1,e/e2}
        \path[edge] (\source) -- (\dest);
    \foreach \source/\dest in {e/g}
        \path[dedge] (\source) -- (\dest);
\node[left] () at (2.5,5) {{\LARGE (a)}};
\end{scope}

\begin{scope}[shift={(11,2)}]
\node[math,rotate=0] () at (1,4) {$\tikzLongrightarrow$};
\foreach \pos/\name in {{(3,6)/b},{(3,4)/c},{(4,3)/f},{(6,3)/e1},{(7,4)/e2}}
        \node[enode] (\name) at \pos {};
   \foreach \source/\dest in {f/c,c/b,b/f,f/e1,e1/e2,c/e1,f/e2}
        \path[edge] (\source) -- (\dest);
\end{scope}

%\node[math,rotate=22.5] () at (6,5) {$\tikzLongrightarrow$};

\begin{scope}[shift={(18,2)}]
  \foreach \pos/\name in {{(3,6)/b},{(3,4)/c},{(4,3)/f},{(6,3.5)/e1},{(6,4.5)/e2},{(7,4)/j}}
        \node[enode] (\name) at \pos {};
  \foreach \pos/\name in {{(4,5)/g},{(5,4)/e}}
        \node[pnode] (\name) at \pos {};
   \foreach \source/\dest in {f/c,g/f,f/e,c/b,c/g,g/b,e/e1,e/e2,j/e1,j/e2,e1/e2}
        \path[edge] (\source) -- (\dest);
    \foreach \source/\dest in {e/g}
        \path[dedge] (\source) -- (\dest);
                \node[left] () at (2.5,5) {{\LARGE (b)}};
%\node[math,right] () at (f.east) {{\LARGE $w$}};
%\node[math,right] () at (g.east) {{\LARGE $u$}};
%\node[math,left] () at (c.west) {{\LARGE $t$}};
%\node[math,left] () at (b.west) {{\LARGE $s$}};
\end{scope}

\begin{scope}[shift={(25,2)}]
\node[math,rotate=0] () at (1,4) {$\tikzLongrightarrow$};
  \foreach \pos/\name in {{(3,6)/b},{(3,4)/c},{(4,3)/f},{(6,3.5)/e1},{(6,4.5)/e2},{(7,4)/j}}
        \node[enode] (\name) at \pos {};
   \foreach \source/\dest in {b/f,c/b,j/e1,j/e2}
        \path[edge] (\source) -- (\dest);
    \foreach \source/\dest in {c/f,e1/e2}
        \path[dedge] (\source) -- (\dest);
\end{scope}

\begin{scope}[shift={(4,-2)}]
\foreach \pos/\name in {{(3,6)/b},{(3,4)/c},{(3,2)/a},{(4,3)/f},{(6,3.5)/e1},{(6,4.5)/e2}}
        \node[enode] (\name) at \pos {};
  \foreach \pos/\name in {{(4,5)/g},{(5,4)/e}}
        \node[pnode] (\name) at \pos {};
   \foreach \source/\dest in {f/c,g/f,f/e,c/b,c/g,g/b,e1/e2,e/e1,e/e2,a/c,a/f}
        \path[edge] (\source) -- (\dest);
    \foreach \source/\dest in {e/g}
        \path[dedge] (\source) -- (\dest);
\node[left] () at (2.5,5) {{\LARGE (c)}};
\end{scope}

\begin{scope}[shift={(11,-2)}]
        \node[math,rotate=0] () at (1,4) {$\tikzLongrightarrow$};
        \foreach \pos/\name in {{(3,6)/b},{(3,4)/c},{(3,2)/a},{(4,3)/f},{(6,3.5)/e1},{(6,4.5)/e2}}
        \node[enode] (\name) at \pos {};
   \foreach \source/\dest in {e1/e2,a/c,a/f,b/e1,b/e2}
        \path[edge] (\source) -- (\dest);
             \foreach \source/\dest in {c/f}
        \path[tredge] (\source) -- (\dest);
\end{scope}

%\node[math,rotate=22.5] () at (6,5) {$\tikzLongrightarrow$};

\begin{scope}[shift={(18,-2)}]
\foreach \pos/\name in {{(3,6)/b},{(3,2)/a},{(4,3)/f},{(6,3.5)/e1},{(6,4.5)/e2}}
        \node[enode] (\name) at \pos {};
  \foreach \pos/\name in {{(4,5)/g},{(3,4)/c},{(5,4)/e}}
        \node[pnode] (\name) at \pos {};
   \foreach \source/\dest in {f/c,g/f,f/e,c/b,c/g,g/b,e1/e2,e/e1,e/e2,a/c,a/f}
        \path[edge] (\source) -- (\dest);
    \foreach \source/\dest in {e/g}
        \path[dedge] (\source) -- (\dest);
        \path[edge,dotted] (b) .. controls +(-.5,-1) .. (c);
                 \node[left] () at (2.5,5) {{\LARGE (d)}};
\end{scope}

\begin{scope}[shift={(25,-2)}]
\node[math,rotate=0] () at (1,4) {$\tikzLongrightarrow$};
\foreach \pos/\name in {{(3,6)/b},{(3,4)/c},{(3,2)/a},{(6,3.5)/e1},{(6,4.5)/e2}}
        \node[enode] (\name) at \pos {};
  \foreach \pos/\name in {{(5,4)/e}}
        \node[pnode] (\name) at \pos {};
   \foreach \source/\dest in {c/b,e1/e2,e/e1,e/e2,a/c,a/e,b/e,c/e}
        \path[edge] (\source) -- (\dest);
        \path[edge,dotted] (b) .. controls +(-.5,-1) .. (c);
 \end{scope}

\begin{scope}[shift={(4,-8)}]
  \foreach \pos/\name in {{(3,6)/b},{(3,4)/c},{(6,3.5)/e1},{(6,4.5)/e2},{(4,1)/g},{(3.5,2)/f1},{(4.5,2)/f2}}
        \node[enode] (\name) at \pos {};
  \foreach \pos/\name in {{(4,5)/h},{(4,3)/f},{(5,4)/e}}
        \node[pnode] (\name) at \pos {};
   \foreach \source/\dest in {f/c,h/f,f/e,c/b,c/h,h/b,e/e1,e/e2,e1/e2,f/f1,f/f2,g/f1,g/f2,f1/f2}
        \path[edge] (\source) -- (\dest);
    \foreach \source/\dest in {e/h}
        \path[dedge] (\source) -- (\dest);
              \node[left] () at (2.5,5) {{\LARGE (e)}};
%\node[math,right] () at (f.east) {{\LARGE $w$}};
%\node[math,right] () at (g.east) {{\LARGE $u$}};
%\node[math,left] () at (c.west) {{\LARGE $t$}};
%\node[math,left] () at (b.west) {{\LARGE $s$}};
\end{scope}

\begin{scope}[shift={(11,-8)}]
\node[math,rotate=0] () at (1,4) {$\tikzLongrightarrow$};
 \foreach \pos/\name in {{(3,6)/b},{(3,4)/c},{(6,3.5)/e1},{(6,4.5)/e2},{(4,1)/g},{(3.5,2)/f1},{(4.5,2)/f2}}
        \node[enode] (\name) at \pos {};
  \foreach \pos/\name in {{(4.5,4)/j}}
        \node[pnode] (\name) at \pos {};
   \foreach \source/\dest in {c/b,j/b,e1/e2,g/f1,g/f2,j/e1,j/e2}
        \path[edge] (\source) -- (\dest);
            \foreach \source/\dest in {c/j,f1/f2}
        \path[dedge] (\source) -- (\dest);
  
   \end{scope}

\begin{scope}[shift={(19,-8)}]
  \foreach \pos/\name in {{(2,5)/a},{(3,6)/b},{(6,3.5)/e1},{(6,4.5)/e2},{(3.5,2)/f1},{(4.5,2)/f2}}
        \node[enode] (\name) at \pos {};
  \foreach \pos/\name in {{(3,4)/c},{(4,5)/h},{(4,3)/f},{(5,4)/e}}
        \node[pnode] (\name) at \pos {};
   \foreach \source/\dest in {a/b,f/c,h/f,f/e,c/b,c/h,h/b,e/e1,e/e2,e1/e2,f/f1,f/f2,f1/f2}
        \path[edge] (\source) -- (\dest);
    \foreach \source/\dest in {e/h,a/c}
        \path[dedge] (\source) -- (\dest);
              \node[left] () at (1.5,5) {{\LARGE (f)}};
\end{scope}

\begin{scope}[shift={(25,-8)}]
\node[math,rotate=0] () at (1,4) {$\tikzLongrightarrow$};
  \foreach \pos/\name in {{(2,5)/a},{(3,6)/b},{(6,3.5)/e1},{(6,4.5)/e2},{(3.5,2)/f1},{(4.5,2)/f2}}
        \node[enode] (\name) at \pos {};
    \foreach \source/\dest in {e1/e2,e2/f2,e1/f2,f1/a,f1/b}
        \path[edge] (\source) -- (\dest);
    \foreach \source/\dest in {a/b}
        \path[dedge] (\source) -- (\dest);
\end{scope}

%\node[math,rotate=0] () at (5,-2) {$\tikzLongrightarrow$};
%
%\begin{scope}[shift={(5,-7)}]
%  \foreach \pos/\name in {{(3,6)/b},{(3,4)/c},{(4,3)/f},{(6,3.5)/e1},{(6,4.5)/e2}}
%        \node[enode] (\name) at \pos {};
%   \foreach \source/\dest in {b/e1,e2/b,e1/e2}
%        \path[edge] (\source) -- (\dest);
%    \foreach \source/\dest in {c/f}
%        \path[tredge] (\source) -- (\dest);
%\end{scope}

\end{tikzpicture}
\end{center}
\caption{One double edge in a diamond}\label{f:onedoub}
\end{figure}
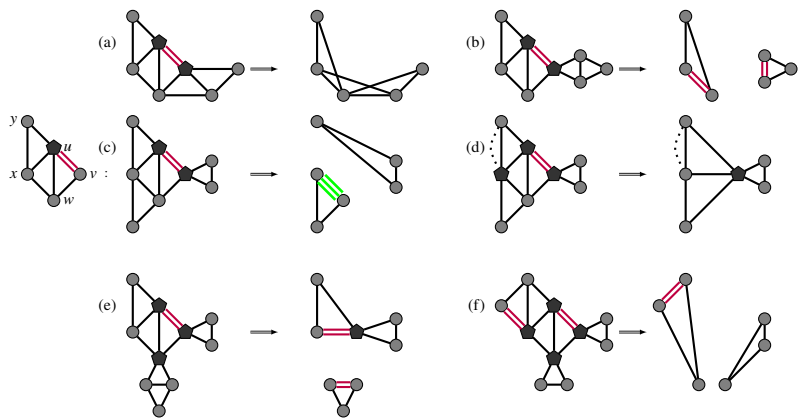

Finally, a double edge can be either in the centre or the outside of a diamond, but in both cases that gives four known edges to at least one vertex in the diamond; if the double edge is in the centre of a diamond then we can use the X-reduction however the fifth edges are attached, using observation \ref{o:four}. Thus we can assume henceforth that no double edge is in the centre of a diamond, and, additionally, we will again use that no edge is in more than two triangles, leading to the configurations shown in figure \ref{f:onedoub}, looking at how the remaining edges can be at $v$.

%
%The top configuration is always reducible in this way, and the middle will work so long as the triangle at $v$ is not aloof. $w$ is not part of a double edge as it is in a diamond with $uv$, and if it is joined to a triangle which is pendant but not aloof then we can use the third reduction. 

If $vw$ is in two triangles then we can use reduction (a), so we can suppose that $v$ has a pendant triangle. If it is not aloof we can use reduction (b), so can now assume that $v$ has an aloof pendant triangle. Similarly, if $wx$ is in two triangles then reduction (c) will work unless edge $xy$ is either part of a second triangle or is a double edge, in which case we use reduction (d), where the dotted arc indicates either option.

Thus we can assume that $w$ has a pendant triangle, and if it isn't aloof then we can use the reduction shown in (e), so now both $v$ and $w$ have aloof pendant triangles.
If $x$ has two neighbours outside of the configuration then we can delete 
$\{u,v,w,x\}$, join the remaining neighbours of $v$ and $w$ in an $A_6$ and then $y$ and $x$'s remaining neighbours can be joined in a triangle.

If $x$ has less than two neighbours outside of the configuration then it must have one, since its only known non-pentagonal neighbour is $y$ and if $xy$ is a triple edge then $y$'s fifth edge cannot be in a triangle. However, $xy$ cannot then be a double edge as it would be in the centre of a diamond, so can be reduced as in (f).

%In this case, if the last two edges at $w$ do not form an aloof triangle $A$, we can remove the edges from $w$ to $A$ and double the third edge of $A$ and then contract $v$, $w$ and $u$ into one vertex and what remains is 5-regular and has the triangle property. Thus both $v$ and $w$ can be assumed to have pendant triangles.
%
%However, we are now able to delete $u$, $v$ and $w$ and, either contract $st$ or delete $t$ and connect its neighbours in triangles to make that section quintic. The four other vertices of pendant triangles that were at $v$ and $w$ can be joined to two new vertices to make $A_6$ as in figure \ref{f:a6red}, thus giving us a reduction overall.

\subsection{No multiple edges, no edges in more than two triangles}

Recall that we can assume that we do not have $K_4$, $W_4$ or $W_5$ in the graph as well.
Any triangle that is not part of a diamond is necessarily atomic, but we can consider the cases of a diamond $H$ with respect to $s_H$ (now more accurately defined as the number of vertices in $G\backslash H$ which are adjacent to two {\em adjacent} vertices of $H$). Using the degrees of $H$, we can see that $0 \leq s_H \leq 5$, and note that no vertex can be adjacent to both vertices in the central edge of $H$, so that actually $s_H \leq 4$, with at most one new vertex adjacent to any pair of adjacent vertices of $H$.

\begin{itemize}
\item{Case i) $0 \leq s_H \leq 1$}\\
For $s_H = 0$ a Z-reduction exists which will necessarily give rise to a smaller quintic graph with the triangle property.
When $s_H=1$ only one configuration is possible since the central edge of the diamond is already in two triangles. However, in this configuration there must be a vertex adjacent to all four other vertices and its fifth edge needs to be in a triangle. This means that we actually have a diamond which has $s_H=2$, so we can move to case iii).
\item{Case ii) $s_H = 2$}\\
If there is no diamond in $G$ which has $s_H > 2$ then there are only two possible configurations as shown in figure \ref{f:sh2}, and, moreover, the remaining edges from the named vertices cannot be part of a triangle with any of the other vertices, and thus the configurations to their right must exist, which can then be reduced as shown.

\begin{figure}[ht]
\begin{center}
\begin{tikzpicture}[scale=0.5, transform shape]
  \foreach \pos/\name in {{(3,8)/a},{(3,6)/b},{(4,7)/c},{(4,5)/d},{(5,6)/e},{(5,4)/f}}
        \node[enode] (\name) at \pos {};
   \foreach \source/\dest in {a/b,a/c,b/c,b/d,c/d,c/e,d/e,d/f,e/f}
        \path[edge] (\source) -- (\dest);
 \node[math,right] () at (c.east) {{\LARGE $v$}};
\node[math,left] () at (d.west) {{\LARGE $w$}};
\node[math,rotate=0] () at (6.5,6) {\Large :};

\begin{scope}[shift={(0,-5)}]
  \foreach \pos/\name in {{(3,4)/a},{(3,6)/b},{(4,7)/c},{(5,6)/e},{(5,4)/f}}
        \node[enode] (\name) at \pos {};
  \foreach \pos/\name in {{(4,5)/d}}
        \node[pnode] (\name) at \pos {};
   \foreach \source/\dest in {a/b,a/d,b/c,b/d,c/d,c/e,d/e,d/f,e/f}
        \path[edge] (\source) -- (\dest);
 \node[math,right] () at (c.east) {{\LARGE $u_1$}};
 \node[math,right] () at (e.east) {{\LARGE $u_2$}};
 \node[math,right] () at (f.east) {{\LARGE $u_3$}};
 \node[math,left] () at (a.west) {{\LARGE $u_4$}};
 \node[math,left] () at (b.west) {{\LARGE $u_0$}};
 \node[math,rotate=0] () at (6.5,6) {\Large :};
\end{scope}

\begin{scope}[shift={(5,0)}]
  \foreach \pos/\name in {{(3,8)/a},{(3,6)/b},{(5,6)/e},{(5,4)/f},{(4,8)/g},{(4,4)/h}}
        \node[enode] (\name) at \pos {};
 \foreach \pos/\name in {{(4,7)/c},{(4,5)/d}}
        \node[pnode] (\name) at \pos {};
  \foreach \source/\dest in {a/b,a/c,b/c,b/d,c/d,c/e,d/e,d/f,e/f,a/g,c/g,h/d,h/f}
        \path[edge] (\source) -- (\dest);
        \node[math,rotate=0] () at (6.5,6) {$\tikzLongrightarrow$};
\end{scope}

\begin{scope}[shift={(5,-5)}]
  \foreach \pos/\name in {{(3,4)/a},{(3,6)/b},{(5,6)/e},{(5,4)/f},{(3.5,8)/g},{(4.5,8)/h}}
        \node[enode] (\name) at \pos {};
  \foreach \pos/\name in {{(4,7)/c},{(4,5)/d}}
        \node[pnode] (\name) at \pos {};
   \foreach \source/\dest in {a/b,a/d,b/c,b/d,c/d,c/e,d/e,d/f,e/f,c/g,c/h,g/h}
        \path[edge] (\source) -- (\dest);
\node[math,rotate=0] () at (6.5,6) {$\tikzLongrightarrow$};
\end{scope}

\begin{scope}[shift={(10,0)}]
  \foreach \pos/\name in {{(3,8)/a},{(3,6)/b},{(5,6)/e},{(5,4)/f},{(4,8)/g},{(4,4)/h}}
        \node[enode] (\name) at \pos {};
  \foreach \source/\dest in {b/g,a/g,h/f,h/e}
        \path[edge] (\source) -- (\dest);
  \foreach \source/\dest in {a/b,e/f}
        \path[dedge] (\source) -- (\dest);
\end{scope}

\begin{scope}[shift={(10,-5)}]
  \foreach \pos/\name in {{(3,4)/a},{(3,6)/b},{(5,6)/e},{(5,4)/f},{(3.5,7)/g},{(4.5,7)/h}}
        \node[enode] (\name) at \pos {};
   \foreach \source/\dest in {a/b,b/e,e/f,f/b,a/e}
        \path[edge] (\source) -- (\dest);
  \foreach \source/\dest in {g/h}
        \path[dedge] (\source) -- (\dest);
\node[math,rotate=0] () at (5.7,6) {\Large or};
\end{scope}

\begin{scope}[shift={(13.5,-5)}]
  \foreach \pos/\name in {{(3,4)/a},{(3,6)/b},{(5,6)/e},{(4,4)/f},{(3.5,7)/g},{(4.5,7)/h}}
        \node[enode] (\name) at \pos {};
   \foreach \source/\dest in {a/b,a/e,f/g,f/h,g/h}
        \path[edge] (\source) -- (\dest);
  \foreach \source/\dest in {b/e}
        \path[dedge] (\source) -- (\dest);
\end{scope}

\end{tikzpicture}
\end{center}
\caption{Reductions when $s_H=2$}\label{f:sh2}
\end{figure}
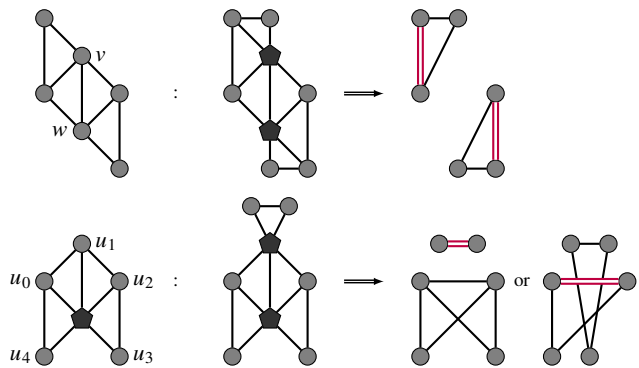

The top reduction will always give a quintic graph with the triangle property, but the first reduction on the bottom will fail if the pendant triangle is aloof and the second if the edge $u_2u_3$ is part of another triangle in $G$. There is a third reduction with the same vertical symmetry as the second, so we can also assume that $u_0u_4$ is part of another triangle in $G$ too. Now we can again use observation \ref{o:four} for $u_0$ and $u_2$ and then reduce as shown in figure \ref{f:sh2b}. The resulting graph will be quintic and will only fail to have the triangle property if $u_3u_4$ is an edge in $G$. In this case we have the 6-wheel $W_6$ and this is another atom; $A_9$.

\begin{figure}[ht]
\begin{center}
\begin{tikzpicture}[scale=0.6, transform shape]

\begin{scope}[shift={(0,0)}]
  \foreach \pos/\name in {{(2,7)/y},{(2,5)/z},{(3,4)/a},{(5,4)/f},{(3.5,8)/g},{(4.5,8)/h},{(6,7)/i},{(6,5)/j}}
        \node[enode] (\name) at \pos {};
  \foreach \pos/\name in {{(4,7)/c},{(4,5)/d},{(3,6)/b},{(5,6)/e}}
        \node[pnode] (\name) at \pos {};
   \foreach \source/\dest in {a/b,a/d,b/c,b/d,c/d,c/e,d/e,d/f,e/f,c/g,c/h,g/h,y/z,y/b,z/a,z/b,e/j,e/i,i/j,f/j}
        \path[edge] (\source) -- (\dest);
\node[math,rotate=0] () at (7.5,6) {$\tikzLongrightarrow$};
\node[math,above] () at (b.north) {{\Large $u_0$}};
\node[math,above] () at (e.north) {{\Large $u_2$}};
\node[math,below] () at (f.south) {{\Large $u_3$}};
\node[math,below] () at (a.south) {{\Large $u_4$}};
\end{scope}

\begin{scope}[shift={(7,0)}]
  \foreach \pos/\name in {{(2,7)/y},{(2,5)/z},{(3,4)/a},{(5,4)/f},{(3.5,8)/g},{(4.5,8)/h},{(6,7)/i},{(6,5)/j}}
        \node[enode] (\name) at \pos {};
   \foreach \source/\dest in {y/z,z/a,a/g,g/z,y/a,f/i,i/j,h/j,h/f,f/j}
        \path[edge] (\source) -- (\dest);
        \node[math,rotate=0] () at (7,6) {\Large unless};
\end{scope}

\begin{scope}[shift={(12.5,.5)},rotate=0]
  \foreach \pos/\name in {{(3,4)/a},{(5,4)/f},{(4,7)/c},{(3,6)/b},{(5,6)/e}}
        \node[enode] (\name) at \pos {};
  \foreach \pos/\name in {{(4,5)/d}}
        \node[pnode] (\name) at \pos {};
   \foreach \source/\dest in {a/f,f/e,e/c,c/b,b/a,d/a,d/b,d/c,d/e,d/f}
        \path[edge] (\source) -- (\dest);
%     \draw \foreach \x in {0,72,...,288}
%         {
%                  (\x:2) -- (\x-72:2)
%                  (\x:2) -- (0:0)
%         };
%\node[enode] () at (0,0) {};
%  \foreach \x in {0,72,...,288} {\node[enode] () at (\x:2) {};}
\end{scope}

\end{tikzpicture}
\end{center}
\caption{Final reduction for $s_H=2$}\label{f:sh2b}
\end{figure}
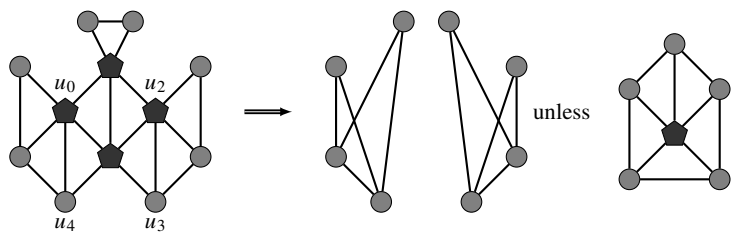

However, each vertex in the rim of $W_6$ must be adjacent to a pendant triangle, using its rotational symmetry, and if any one of these pendant triangles is not aloof we can delete the whole $W_6$, and use lemma \ref{l:2s} to double the edges from the non-aloof pendant triangles and rejoin the vertices of the aloof triangles to 2 or 4 new vertices of degree 5. Again, atom $A_9$ is only not reducible in this way if all of its pendant triangles are aloof.
However, we can instead use two copies of $A_3$ together with their two pendant triangles to cover the ten vertices of degree 2, and this will give a smaller quintic graph with the triangle property, so $A_9$ is reducible.

\item{Case iii) $s_H = 3$}\\
By symmetry, there is only one way to have $s_H=3$. Taking $H$ to be the subgraph induced by $\{t,u,v,w\}$, vertex $v$ in figure \ref{f:shb3} must be in a new triangle as shown, otherwise $s_H=4$. We can use a Z-reduction using $vuwx$ unless $xw$ is in another triangle, but in that case $u$'s fifth neighbour has to be a new vertex as in the configuration on the right and thus we can reduce using double edges as shown.

\begin{figure}[ht]
\begin{center}
\begin{tikzpicture}[scale=0.5, transform shape]

\begin{scope}[shift={(0,0)}]
  \foreach \pos/\name in {{(3,8)/a},{(4,9)/b},{(5,10)/c},{(6,9)/f},{(5,6)/g},{(6,7)/h}}
        \node[enode] (\name) at \pos {};
 \foreach \pos/\name in {{(4,7)/d},{(5,8)/e}}
        \node[pnode] (\name) at \pos {};
   \foreach \source/\dest in {a/b,b/c,a/d,b/d,b/e,c/e,c/f,f/e,d/e,d/g,d/h,e/h,g/h}
        \path[edge] (\source) -- (\dest);
 \node[math,right] () at (b.east) {{\LARGE $t$}};
 \node[math,right] () at (e.east) {{\LARGE $v$}};
\node[math,right] () at (h.east) {{\LARGE $w$}};
\node[math,left] () at (d.west) {{\LARGE $u$}};
\node[math,left] () at (g.west) {{\LARGE $x$}};
\node[math,rotate=0] () at (7,8) {$\tikzLongrightarrow$};
\end{scope}

\begin{scope}[shift={(5.5,0)}]
  \foreach \pos/\name in {{(3,8)/a},{(4,9)/b},{(5,10)/c},{(5,6)/d},{(6,7)/e},{(6,9)/f}}
        \node[enode] (\name) at \pos {};
   \foreach \source/\dest in {a/b,b/c,a/d,b/d,b/e,c/e,c/f,f/e}
        \path[edge] (\source) -- (\dest);
\end{scope}

\begin{scope}[shift={(10,0)}]
  \foreach \pos/\name in {{(3,8)/a},{(4,9)/b},{(5,10)/c},{(6,9)/f},{(5,6)/g},{(6,5)/i},{(7,6)/j}}
        \node[enode] (\name) at \pos {};
 \foreach \pos/\name in {{(4,7)/d},{(5,8)/e},{(6,7)/h}}
        \node[pnode] (\name) at \pos {};
   \foreach \source/\dest in {a/b,b/c,a/d,b/d,b/e,c/e,c/f,f/e,d/e,d/g,d/h,e/h,g/i,h/i,h/j,i/j,g/h}
        \path[edge] (\source) -- (\dest);
\node[math,rotate=0] () at (7,8) {$\tikzLongrightarrow$};
\end{scope}

\begin{scope}[shift={(15.5,0)}]
  \foreach \pos/\name in {{(3,8)/a},{(4,9)/b},{(5,10)/c},{(6,9)/f},{(4,6)/g},{(5,5)/i},{(6,6)/j}}
        \node[enode] (\name) at \pos {};
 \foreach \pos/\name in {{(5,8)/e}}
        \node[pnode] (\name) at \pos {};
   \foreach \source/\dest in {a/b,b/c,c/e,c/f,f/e,i/j,a/e,g/j}
        \path[edge] (\source) -- (\dest);
  \foreach \source/\dest in {b/e,g/i}
        \path[dedge] (\source) -- (\dest);
\end{scope}

\end{tikzpicture}
\end{center}
\caption{Reductions when $s_H=3$}\label{f:shb3}
\end{figure}
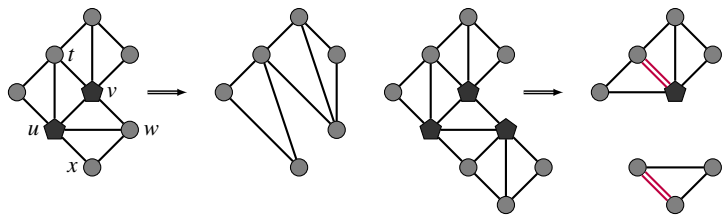

\item{Case iv) $s_H = 4$}\\
Finally, we can assume that $H$ has four new vertices, one adjacent to each of its edges, as shown in figure \ref{f:sh4}. Vertices $v$ and $w$ need to have their fifth edge adjacent to one of the neighbours of degree 2 in the configuration and this gives a configuration such as is in the figure with degree sequence $\{5,5,5,5,3,3,2,2,2,2\}$ with the only edges between the vertices not of degree 5 being between each 3 and one 2, however this is done. This can be reduced to the configuration shown which is certainly quintic and has the triangle property.

\begin{figure}[ht]
\begin{center}
\begin{tikzpicture}[scale=0.5, transform shape]
\begin{scope}[shift={(0,0)}]
  \foreach \pos/\name in {{(3,6)/a},{(3,8)/b},{(3,10)/c},{(5,6)/f},{(5,8)/g},{(5,10)/h}}
        \node[enode] (\name) at \pos {};
 \foreach \pos/\name in {{(4,7)/d},{(4,9)/e}}
        \node[pnode] (\name) at \pos {};
   \foreach \source/\dest in {a/b,b/c,a/d,b/d,b/e,c/e,d/e,h/e,d/e,h/g,g/e,g/f,d/f,d/g}
        \path[edge] (\source) -- (\dest);
 \node[math,right] () at (g.east) {{\LARGE $w$}};
\node[math,left] () at (b.west) {{\LARGE $v$}};

\node[math,rotate=0] () at (6,8) {\LARGE :};
\end{scope}

\begin{scope}[shift={(5,0)}]
  \foreach \pos/\name in {{(3,6)/a},{(3,10)/c},{(5,6)/f},{(5,10)/h},{(2,9)/i},{(6,9)/j}}
        \node[enode] (\name) at \pos {};
 \foreach \pos/\name in {{(3,8)/b},{(4,7)/d},{(4,9)/e},{(5,8)/g}}
        \node[pnode] (\name) at \pos {};
   \foreach \source/\dest in {a/b,b/c,a/d,b/d,b/e,c/e,d/e,h/e,d/e,h/g,g/e,g/f,d/f,d/g,b/i,c/i,g/j,h/j}
        \path[edge] (\source) -- (\dest);
\node[math,rotate=0] () at (7,8) {$\tikzLongrightarrow$};
\end{scope}

\begin{scope}[shift={(10,0)}]
  \foreach \pos/\name in {{(3,6)/a},{(3,10)/c},{(5,6)/f},{(5,10)/h},{(2,9)/i},{(6,9)/j}}
        \node[enode] (\name) at \pos {};
 \foreach \pos/\name in {{(4,7)/d},{(4,9)/e}}
        \node[pnode] (\name) at \pos {};
   \foreach \source/\dest in {a/d,c/e,d/e,h/e,d/e,d/f,c/i,h/j,i/e,j/e,c/d,h/d,a/f}
        \path[edge] (\source) -- (\dest);
\end{scope}

\end{tikzpicture}
\end{center}
\caption{Reduction for $s_H=4$}\label{f:sh4}
\end{figure}
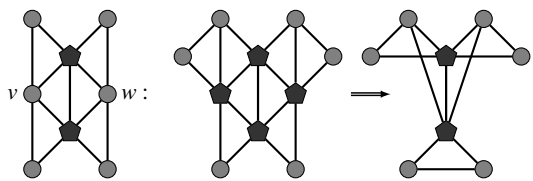
\end{itemize}

\section{Atomic construction}\label{s:atom}

Throughout the previous work we have discovered various atoms which were not reducible using the basic methods introduced. 
%These included $K_3$ and $K_4$ which were the base for the fundamental graphs in subsection \ref{s:irred}. 
Each of the remaining atoms contains a combination of possibly some degree 5 vertices along with other vertices of degree 2 (taking into account the proven necessary pendant aloof triangles) or degree 3 as shown in table \ref{t:atoms23}. 
We have already removed $A_5$, $A_6$, $A_8$ and $A_9$ from the table since they were shown reducible in earlier sections.% \ref{s:reduc}.

\begin{table}[ht]
\begin{center}
\begin{tabular}{|c|c|c|c|c|}\hline
Atom& Configuration & Degree 5& Degree 3 & Degree 2 \\\hline
$A_1$ & \begin{tikzpicture}[scale=0.3, transform shape]
  \foreach \pos/\name in {{(0,0)/a},{(2,0)/b},{(1,1.7)/c}}
        \node[enode] (\name) at \pos {};
 \foreach \pos/\name in {}
        \node[pnode] (\name) at \pos {};
   \foreach \source/\dest in {a/b,b/c,a/c}
        \path[edge] (\source) -- (\dest);
\end{tikzpicture} & 0 & 0 & 3 \\\hline
$A_2$ & \begin{tikzpicture}[scale=0.3, transform shape]
  \foreach \pos/\name in {{(0,0)/a},{(2,0)/b},{(0,2)/c},{(2,2)/d}}
        \node[enode] (\name) at \pos {};
 \foreach \pos/\name in {}
        \node[pnode] (\name) at \pos {};
   \foreach \source/\dest in {a/b,b/c,a/c,a/d,b/d,c/d}
        \path[edge] (\source) -- (\dest);
\end{tikzpicture} & 0 & 4 & 0  \\\hline
$A_3$ &  \begin{tikzpicture}[scale=0.3, transform shape]
 \foreach \pos/\name in {{(0,0)/a},{(2,0)/b},{(1,1.7)/c}}
        \node[enode] (\name) at \pos {};
   \foreach \source/\dest in {b/c,a/c}
        \path[edge] (\source) -- (\dest);
          \foreach \source/\dest in {a/b}
        \path[dedge] (\source) -- (\dest);
        \end{tikzpicture} & 0 & 2 & 1 \\\hline
   %larger $A_3$ &  \begin{tikzpicture}[scale=0.3, transform shape]
%\foreach \pos/\name in {{(-2,0)/d},{(4,0)/e},{(1,1.7)/c},{(-1,-1.7)/f},{(3,-1.7)/g}}
%        \node[enode] (\name) at \pos {};
% \foreach \pos/\name in {{(0,0)/a},{(2,0)/b}}
%        \node[pnode] (\name) at \pos {};
%   \foreach \source/\dest in {b/c,a/c,a/d,a/f,f/d,b/e,b/g,e/g}
%        \path[edge] (\source) -- (\dest);
%          \foreach \source/\dest in {a/b}
%        \path[dedge] (\source) -- (\dest);
%        \end{tikzpicture} & 2 & 0 & 5 \\\hline
        $A_4$ &  \begin{tikzpicture}[scale=0.3, transform shape]
         \foreach \pos/\name in {{(1,1.7)/c}}
        \node[enode] (\name) at \pos {};
 \foreach \pos/\name in {{(0,0)/a},{(2,0)/b}}
        \node[pnode] (\name) at \pos {};
   \foreach \source/\dest in {b/c,a/c}
        \path[edge] (\source) -- (\dest);  
        \foreach \source/\dest in {b/a}
        \path[quple] (\source) -- (\dest);
\end{tikzpicture} & 2 & 0 & 1 \\\hline
%        $A_6$ &  \begin{tikzpicture}[scale=0.3, transform shape]
%          \foreach \pos/\name in {{(0,0)/a},{(4,0)/c},{(0,2)/d},{(4,2)/f}}
%        \node[enode] (\name) at \pos {};
% \foreach \pos/\name in {{(2,2)/e},{(2,0)/b}}
%        \node[pnode] (\name) at \pos {};
%   \foreach \source/\dest in {b/e,a/b,b/c,a/e,c/e,d/b,d/e,f/b,f/e}
%        \path[edge] (\source) -- (\dest);
%\end{tikzpicture} & 2 & 0 & 4 \\\hline
        $A_7$ &  \begin{tikzpicture}[scale=0.3, rotate=90, transform shape]
          \foreach \pos/\name in {{(1,1.7)/c},{(1,-1.7)/d}}
        \node[enode] (\name) at \pos {};
 \foreach \pos/\name in {{(0,0)/a},{(2,0)/b}}
        \node[pnode] (\name) at \pos {};
   \foreach \source/\dest in {b/c,a/c,a/d,b/d}
        \path[edge] (\source) -- (\dest);
         \foreach \source/\dest in {a/b}
        \path[tredge] (\source) -- (\dest);
\end{tikzpicture} & 2 & 0 & 2 \\\hline
     $A_{10}$ &  \begin{tikzpicture}[scale=0.3, transform shape]
      \foreach \pos/\name in {{(0,0)/a}}
        \node[enode] (\name) at \pos {};
 \foreach \pos/\name in {{(0,2)/c},{(2,0)/b},{(2,2)/d}}
        \node[pnode] (\name) at \pos {};
   \foreach \source/\dest in {a/b,a/c,a/d}
        \path[edge] (\source) -- (\dest);
               \foreach \source/\dest in {c/d,d/b,c/b}
        \path[dedge] (\source) -- (\dest);
\end{tikzpicture} & 3 & 1 & 0 \\\hline
     $A_{11}$ &  \begin{tikzpicture}[scale=0.3, transform shape]
      \foreach \pos/\name in {{(0,0)/a}}
        \node[enode] (\name) at \pos {};
 \foreach \pos/\name in {{(0,2)/c},{(2,0)/b},{(2,2)/d}}
        \node[pnode] (\name) at \pos {};
   \foreach \source/\dest in {a/c,c/b}
        \path[edge] (\source) -- (\dest);
               \foreach \source/\dest in {a/b,d/b}
        \path[dedge] (\source) -- (\dest);
                    \foreach \source/\dest in {c/d}
        \path[tredge] (\source) -- (\dest);
\end{tikzpicture} & 3 & 1& 0  \\\hline
\end{tabular}
\end{center}
\caption{Atom list with degree 2 and 3 vertices}\label{t:atoms23}
\end{table}

\subsection{Reductions for cut vertices}\label{s:cutred}
Let $A_{10}$ and $A_{11}$ be the configurations with cut vertices from figure \ref{f:cutv} along with $A_4$. 
When combined together they gave the irreducible graphs shown in figure \ref{f:ircut}.

If $G$ is a graph of connectivity 1 which has $v$ as a cutvertex and neither component of $G-v$ is from $C:=\{A_4, A_{10}, A_{11} \}$ then we can form two smaller quintic graphs with the triangle property by using configurations from $C$ to replace
components of $G-v$. Thus it just remains to show that any graph with a pendant configuration from $C$ can be reduced:

\begin{lem}\label{l:cut3}
If $G$ is a quintic graph with the triangle property that has $A_{10}$ or $A_{11}$ pendant at a cut vertex then $G$ is reducible.
\end{lem}

\begin{proof}
Suppose $G$ is such a graph and let $J$ be $G$ with all four vertices of the pendant atom removed. Necessarily $J$ will have two vertices of degree 4 which must be joined by an edge $e$ because $G$ has the triangle property. If $e$ is in more than one triangle in $G$ it will be in at least one in $J$ and so we can add a multiple edge parallel to $e$ and the resulting graph will be quintic and have the triangle property.

If $e$ is in only one triangle then we can delete $e$ and add a copy of $A_7$ in its place. This graph will have two fewer vertices than $G$ but still be quintic and have the triangle property. If $e$ is a multiple edge we can just remove one of its edges and add a copy of $A_7$ which still satisfies our properties.
\end{proof}

The case for the remaining type of cut vertex is slightly more involved:

\begin{thm}\label{t:a4red}
If $G$ is a quintic graph with the triangle property that has $A_{4}$ pendant at a cut vertex then $G$ is reducible.
\end{thm}

\begin{proof}
Suppose $G$ is such a graph and let $u$ be the cut-vertex and suppose $L:=G-A_4$. $L$ will now not necessarily have the triangle property and will have either two joined vertices of degree 4 and 3 (if there was a double edge at $v$, in a similar way to lemma \ref{l:cut3}) or three vertices of degree 4 joined in a path (otherwise).

\begin{itemize}
\item{Case i) $L$ has two vertices of degree less than 5}\\
Let $v$ be the vertex of $L$ of degree 3 and $w$ the vertex of degree 4.
In $G$ there can be one or two neighbours of $v$ other than $u$ and $w$, there cannot be zero since $u$'s neighbours are known and if $vw$ was a triple edge then the fifth edge from $w$ cannot be part of a triangle. If $v$ and $w$ have no common neighbour then in $L$ we can contract $vw$ and the resulting graph is quintic and has the triangle property.

If $v$ and $w$ have two common neighbours other than $u$ then we can reduce as in a) of figure \ref{f:L2}, and similarly as in b) if $vw$ is a double edge. Otherwise, we have the situations in c) or d) which can be reduced deleting all pentagonal vertices and by inserting an $A_7$ respectively.

\begin{figure}[ht]
\begin{center}
\begin{tikzpicture}[scale=0.45, transform shape]
 \begin{scope}[shift={(0,0)}]
 \node[math,left] () at (2,0) {{\LARGE a)}};
   \foreach \pos/\name in {{(6,0)/c},{(5.5,-1)/d},{(7,0)/e}}
       \node[enode] (\name) at \pos {};
     \foreach \pos/\name in {{(4.5,0)/a},{(5.5,1)/b}, {(3,.8)/y}, {(3,-.8)/z}}
        \node[pnode] (\name) at \pos {};
     \foreach \source/\dest in {a/y,a/z}
        \path[edge] (\source) -- (\dest);
\foreach \source/\dest in {y/z}
        \path[quple] (\source) -- (\dest);
         \foreach \source/\dest in {a/d,b/d,b/c,d/c,b/e,d/e}
        \path[edge] (\source) -- (\dest);
\foreach \source/\dest in {a/b}
        \path[dedge] (\source) -- (\dest);
  \node[math,right] () at (a.east) {{\LARGE $u$}};
   \node[math,right] () at (b.east) {{\LARGE $v$}};
 \node[math,right] () at (d.east) {{\LARGE $w$}};
        \end{scope}
         \begin{scope}[shift={(6,0)}]
 \node[math,left] () at (3,0) {{$\tikzLongrightarrow$}};
   \foreach \pos/\name in {{(6,0)/c},{(5.5,-1)/d},{(7,0)/e}}
       \node[enode] (\name) at \pos {};
     \foreach \pos/\name in {{(4,.8)/y}, {(4,-.8)/z}}
        \node[pnode] (\name) at \pos {};
    \foreach \source/\dest in {d/y,d/z,d/e,c/e}
        \path[edge] (\source) -- (\dest);
\foreach \source/\dest in {y/z}
        \path[quple] (\source) -- (\dest);
         \foreach \source/\dest in {d/c}
        \path[edge] (\source) -- (\dest);
        \end{scope}
        
         \begin{scope}[shift={(13,0)}]
 \node[math,left] () at (2,0) {{\LARGE b)}};
   \foreach \pos/\name in {{(6.5,0)/c},{(5.5,-1)/d}}
       \node[enode] (\name) at \pos {};
     \foreach \pos/\name in {{(4.5,0)/a},{(5.5,1)/b}, {(3,.8)/y}, {(3,-.8)/z}}
        \node[pnode] (\name) at \pos {};
     \foreach \source/\dest in {a/y,a/z}
        \path[edge] (\source) -- (\dest);
\foreach \source/\dest in {y/z}
        \path[quple] (\source) -- (\dest);
  \foreach \source/\dest in {a/d,b/c,d/c}
        \path[edge] (\source) -- (\dest);
\foreach \source/\dest in {a/b,d/b}
        \path[dedge] (\source) -- (\dest);
        \end{scope}
         \begin{scope}[shift={(18,0)}]
 \node[math,left] () at (3.5,0) {{$\tikzLongrightarrow$}};
   \foreach \pos/\name in {{(6.5,0)/c},{(5.5,-1)/d}}
       \node[enode] (\name) at \pos {};
     \foreach \pos/\name in {{(4.5,0)/a},{(5.5,1)/b}}
        \node[pnode] (\name) at \pos {};
          \foreach \source/\dest in {b/d,b/c,d/c}
        \path[edge] (\source) -- (\dest);
\foreach \source/\dest in {a/d}
        \path[dedge] (\source) -- (\dest);
 \foreach \source/\dest in {a/b}
        \path[tredge] (\source) -- (\dest);
        \end{scope}
        
         \begin{scope}[shift={(0,-4)}]
 \node[math,left] () at (2,0) {{\LARGE c)}};
   \foreach \pos/\name in { {(6.5,0)/c},{(5.5,-1)/d},{(6.5,2)/e}}
       \node[enode] (\name) at \pos {};
     \foreach \pos/\name in {{(4.5,0)/a},{(5.5,1)/b},{(3,.8)/y}, {(3,-.8)/z}}
        \node[pnode] (\name) at \pos {};
     \foreach \source/\dest in {a/y,a/z}
        \path[edge] (\source) -- (\dest);
\foreach \source/\dest in {y/z}
        \path[quple] (\source) -- (\dest);
         \foreach \source/\dest in {a/d,b/d,b/c,d/c,b/e,c/e}
        \path[edge] (\source) -- (\dest);
\foreach \source/\dest in {a/b}
        \path[dedge] (\source) -- (\dest);
        \end{scope}
         \begin{scope}[shift={(5,-4)}]
 \node[math,left] () at (4,0) {{$\tikzLongrightarrow$}};
    \foreach \pos/\name in { {(6.5,0)/c},{(5.5,-1)/d},{(6.5,2)/e}}
       \node[enode] (\name) at \pos {};
     \foreach \source/\dest in {c/e,d/e}
        \path[edge] (\source) -- (\dest);
\foreach \source/\dest in {c/d}
        \path[dedge] (\source) -- (\dest);
        \end{scope}
        
         \begin{scope}[shift={(13,-4)}]
 \node[math,left] () at (2,0) {{\LARGE d)}};
   \foreach \pos/\name in { {(6.5,0)/c},{(5.5,-1)/d}}
       \node[enode] (\name) at \pos {};
     \foreach \pos/\name in {{(4.5,0)/a},{(5.5,1)/b},{(3,.8)/y}, {(3,-.8)/z}}
        \node[pnode] (\name) at \pos {};
     \foreach \source/\dest in {a/y,a/z}
        \path[edge] (\source) -- (\dest);
\foreach \source/\dest in {y/z}
        \path[quple] (\source) -- (\dest);
         \foreach \source/\dest in {a/d,b/d,d/c}
        \path[edge] (\source) -- (\dest);
\foreach \source/\dest in {a/b,b/c}
        \path[dedge] (\source) -- (\dest);
        \end{scope}
         \begin{scope}[shift={(18,-4)}]
 \node[math,left] () at (4,0) {{$\tikzLongrightarrow$}};
   \foreach \pos/\name in { {(6.5,0)/c},{(5.5,-1)/d}}
       \node[enode] (\name) at \pos {};
     \foreach \pos/\name in {{(4.5,0)/a},{(5.5,1)/b}}
        \node[pnode] (\name) at \pos {};
    \foreach \source/\dest in {a/d,b/d,b/c,d/c,a/c}
        \path[edge] (\source) -- (\dest);
\foreach \source/\dest in {a/b}
        \path[tredge] (\source) -- (\dest);
        \end{scope}
\end{tikzpicture}
\end{center}
\caption{Reductions for a cut vertex with a double edge}\label{f:L2}
\end{figure}
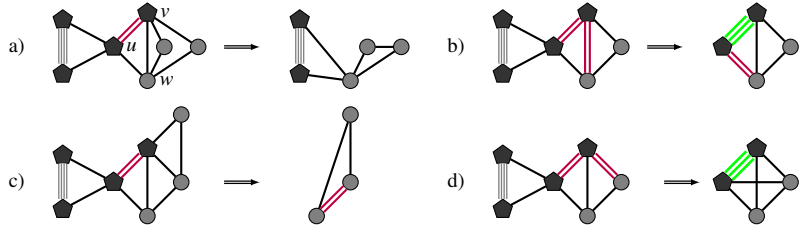
 
\item{Case ii) $L$ has three vertices of degree less than 5}\\
Suppose that $y$ is the vertex of degree 4 in $L$ adjacent to the other two vertices of degree 4, which are $x$ and $z$. 
If $xz$ is an edge then $x$, $y$ and $z$ are indistinguishable and we can reduce as in a) in figure \ref{f:L3} and it will have the triangle property unless each pair of $x$, $y$ and $z$ have a distinct neighbour outside of the configuration as we can use these three vertices in any order. However, in that case we can either reduce as in b), or even remove all pentagonal vertices and insert a triangle.

\begin{figure}[ht]
\begin{center}
\begin{tikzpicture}[scale=0.45, transform shape]
        \begin{scope}[shift={(0,4.5)}]
 \node[math,left] () at (2,0) {{\LARGE a)}};
   \foreach \pos/\name in {{(6.5,0)/c},{(5.5,1)/b},{(5.5,-1)/d}}
       \node[enode] (\name) at \pos {};
     \foreach \pos/\name in {{(4.5,0)/a}, {(3,.8)/y}, {(3,-.8)/z}}
        \node[pnode] (\name) at \pos {};
     \foreach \source/\dest in {a/y,a/z}
        \path[edge] (\source) -- (\dest);
\foreach \source/\dest in {y/z}
        \path[quple] (\source) -- (\dest);
         \foreach \source/\dest in {a/d,b/d,b/c,d/c,a/c,a/b}
        \path[edge] (\source) -- (\dest);
         \node[math,right] () at (b.east) {{\LARGE $x$}};
 \node[math,right] () at (c.east) {{\LARGE $y$}};
 \node[math,right] () at (d.east) {{\LARGE $z$}};
 \end{scope}
         \begin{scope}[shift={(5,4.5)}]
 \node[math,left] () at (4,0) {{$\tikzLongrightarrow$}};
   \foreach \pos/\name in {{(6.5,0)/c},{(5.5,1)/b},{(5.5,-1)/d}}
       \node[enode] (\name) at \pos {};
     \foreach \pos/\name in {{(4.5,0)/a}}
        \node[pnode] (\name) at \pos {};
\foreach \source/\dest in {a/b,a/d}
        \path[dedge] (\source) -- (\dest);
         \foreach \source/\dest in {b/c,d/c,a/c}
        \path[edge] (\source) -- (\dest);
        \end{scope}
           \begin{scope}[shift={(13,4.5)}]
 \node[math,left] () at (2,0) {{\LARGE b)}};
   \foreach \pos/\name in {{(7.3,0)/g},{(6.7,1.7)/e},{(6.7,-1.7)/f}}
       \node[enode] (\name) at \pos {};
     \foreach \pos/\name in {{(4.5,0)/a}, {(3,.8)/y}, {(3,-.8)/z},{(6.2,0)/c},{(5.5,1)/b},{(5.5,-1)/d}}
        \node[pnode] (\name) at \pos {};
     \foreach \source/\dest in {a/y,a/z}
        \path[edge] (\source) -- (\dest);
\foreach \source/\dest in {y/z}
        \path[quple] (\source) -- (\dest);
         \foreach \source/\dest in {a/d,b/d,b/c,d/c,a/c,a/b,b/e,c/e,c/f,d/f,b/g,d/g}
        \path[edge] (\source) -- (\dest);
 \end{scope}
         \begin{scope}[shift={(17,4.5)}]
 \node[math,left] () at (5,0) {{$\tikzLongrightarrow$}};
   \foreach \pos/\name in {{(7.3,0)/g},{(6.7,1.7)/e},{(6.7,-1.7)/f}}
       \node[enode] (\name) at \pos {};
     \foreach \pos/\name in {{(5.5,1)/b},{(5.5,-1)/d}}
        \node[pnode] (\name) at \pos {};
         \foreach \source/\dest in {b/e,d/e,b/f,d/f,b/g,d/g}
        \path[edge] (\source) -- (\dest);
    \foreach \source/\dest in {b/d}
        \path[dedge] (\source) -- (\dest);
        \end{scope}     
 \begin{scope}[shift={(0,0)}]
 \node[math,left] () at (2,0) {{\LARGE c)}};
   \foreach \pos/\name in {{(6.5,0)/c},{(5.5,1)/b},{(5.5,-1)/d}}
       \node[enode] (\name) at \pos {};
     \foreach \pos/\name in {{(4.5,0)/a}, {(3,.8)/y}, {(3,-.8)/z}}
        \node[pnode] (\name) at \pos {};
     \foreach \source/\dest in {a/y,a/z}
        \path[edge] (\source) -- (\dest);
\foreach \source/\dest in {y/z}
        \path[quple] (\source) -- (\dest);
         \foreach \source/\dest in {a/d,b/c,d/c,a/c,a/b}
        \path[edge] (\source) -- (\dest);

  \end{scope}
         \begin{scope}[shift={(6,0)}]
 \node[math,left] () at (3,0) {{$\tikzLongrightarrow$}};
   \foreach \pos/\name in {{(5.5,.5)/b},{(5.5,-1)/d}}
       \node[enode] (\name) at \pos {};
        \foreach \pos/\name in  {{(4,.8)/y}, {(4,-.8)/z}}
        \node[pnode] (\name) at \pos {};
     \foreach \source/\dest in {d/y,d/z}
        \path[edge] (\source) -- (\dest);
        \foreach \source/\dest in {y/z}
        \path[quple] (\source) -- (\dest);
         \node[math,right] () at (b.east) {{\LARGE $\{x,y\}$}};
          \node[math,right] () at (d.east) {{\LARGE $z$}};
       \end{scope}
         \begin{scope}[shift={(13,0)}]
 \node[math,left] () at (2,0) {{\LARGE d)}};
    \foreach \pos/\name in {{(5.5,1)/b},{(5.5,-1)/d}}
       \node[enode] (\name) at \pos {};
     \foreach \pos/\name in {{(4.5,0)/a},{(6.5,0)/c},{(3,.8)/y}, {(3,-.8)/z}}
        \node[pnode] (\name) at \pos {};
     \foreach \source/\dest in {a/y,a/z}
        \path[edge] (\source) -- (\dest);
\foreach \source/\dest in {y/z}
        \path[quple] (\source) -- (\dest);
         \foreach \source/\dest in {a/d,a/b,a/c}
        \path[edge] (\source) -- (\dest);
              \foreach \source/\dest in {b/c,c/d}
        \path[dedge] (\source) -- (\dest);
        \end{scope}
         \begin{scope}[shift={(19,0)}]
 \node[math,left] () at (2.5,0) {{$\tikzLongrightarrow$}};
   \foreach \pos/\name in {{(5.5,1)/b},{(5.5,-1)/d}}
       \node[enode] (\name) at \pos {};
     \foreach \pos/\name in {{(3,.8)/y}, {(3,-.8)/z}}
        \node[pnode] (\name) at \pos {};
 \foreach \source/\dest in {y/z}
        \path[tredge] (\source) -- (\dest);
    \foreach \source/\dest in {b/d,b/y,b/z,d/y,d/z}
        \path[edge] (\source) -- (\dest);
        \end{scope}
        
         \begin{scope}[shift={(0,-4.5)}]
 \node[math,left] () at (2,0) {{\LARGE e)}};
   \foreach \pos/\name in {{(5.5,1)/b},{(5.5,-1)/d},{(6.7,1.7)/e}}
       \node[enode] (\name) at \pos {};
     \foreach \pos/\name in {{(4.5,0)/a},{(6.5,0)/c},{(3,.8)/y}, {(3,-.8)/z}}
        \node[pnode] (\name) at \pos {};
     \foreach \source/\dest in {a/y,a/z}
        \path[edge] (\source) -- (\dest);
\foreach \source/\dest in {y/z}
        \path[quple] (\source) -- (\dest);
         \foreach \source/\dest in {a/d,b/c,a/b,a/c,b/e,c/e}
        \path[edge] (\source) -- (\dest);
              \foreach \source/\dest in {c/d}
        \path[dedge] (\source) -- (\dest);
        \end{scope}
         \begin{scope}[shift={(5.5,-4.5)}]
 \node[math,left] () at (3.5,0) {{$\tikzLongrightarrow$}};
   \foreach \pos/\name in {{(4.5,1)/b},{(4.5,-1)/d},{(5.7,1.7)/e}}
       \node[enode] (\name) at \pos {};
   \foreach \source/\dest in {b/e,d/e}
        \path[edge] (\source) -- (\dest);
   \foreach \source/\dest in {b/d}
        \path[dedge] (\source) -- (\dest);
        \end{scope}
        
         \begin{scope}[shift={(13,-4.5)}]
 \node[math,left] () at (2,0) {{\LARGE f)}};
 \foreach \pos/\name in {{(5.5,1)/b},{(5.5,-1)/d},{(6.7,1.7)/e},{(6.7,-1.7)/f}}
       \node[enode] (\name) at \pos {};
     \foreach \pos/\name in {{(4.5,0)/a},{(6.5,0)/c},{(3,.8)/y}, {(3,-.8)/z}}
        \node[pnode] (\name) at \pos {};
     \foreach \source/\dest in {a/y,a/z}
        \path[edge] (\source) -- (\dest);
\foreach \source/\dest in {y/z}
        \path[quple] (\source) -- (\dest);
         \foreach \source/\dest in {a/d,c/d,b/c,a/b,a/c,b/e,c/e,d/f,c/f}
        \path[edge] (\source) -- (\dest);
                \end{scope}
         \begin{scope}[shift={(17.5,-4.5)}]
 \node[math,left] () at (4,0) {{$\tikzLongrightarrow$}};
 \foreach \pos/\name in {{(4.5,1)/b},{(4.5,-1)/d},{(6.7,1.7)/e},{(6.7,-1.7)/f}}
       \node[enode] (\name) at \pos {};
         \foreach \source/\dest in {b/d,e/d,b/e,b/f,d/f}
        \path[edge] (\source) -- (\dest);
      \end{scope}

\end{tikzpicture}
\end{center}
\caption{Reductions for a cut vertex without a double edge}\label{f:L3}
\end{figure}
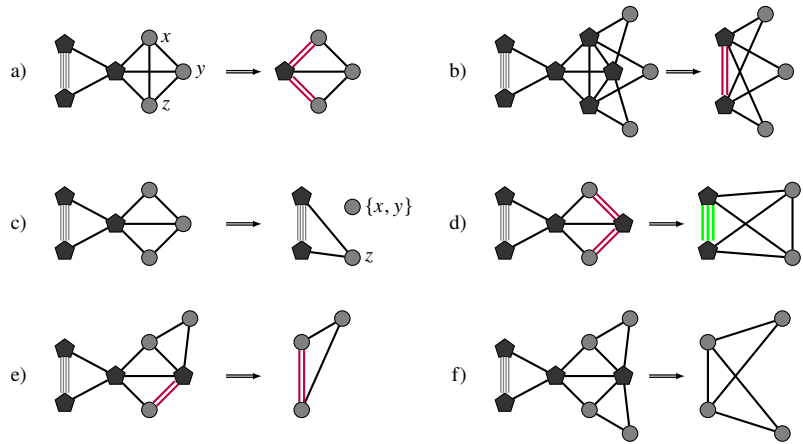

Now we can suppose that $xz$ is not an edge; if $y$ has no multiple edges or  common neighbours with $x$ and $z$ outside the configuration then we can use the Z-reduction as in c) in figure \ref{f:L3}.
Otherwise both $x$ and $z$ either have multiple edges to $y$ or are in a triangle with $y$ outside of the configuration.
%and it doesn't matter if $xz$ is an edge or not.
 
We can simplify by using the symmetry between $x$ and $z$. A triple edge from $x$, say, to $y$ is not possible since the fifth edge from $x$ cannot be in a triangle.
If $xy$ and $yz$ are double edges we can use $A_5$ as shown in d), if only one of these is double then all pentagonal vertices can be deleted and edges added as in e) and if neither are double then we can use $A_7$ as in f).
\end{itemize}
\end{proof}
%We can create all graphs with the triangle property that have just one vertex of degree 2 and all others with degree 5 as arising from a quintic graph $G$ with the triangle property.
%Since the two neighbours of the vertex of degree 2 are either adjacent or not they arise from a double edge in a triangle in $G$ which is subdivided or they come from replacing an $A_7$ subgraph in $G$ by a triangle. Thus we can reduce both $A_{10}$ and $A_{11}$ this way.

\subsection{Final simplifications}\label{s:final}

%Recall that $A_2$ can be reduced using $A_1$ and $A_3$ because $A_2$ has to be surrounded by aloof triangles in order to be not immediately reducible.
%We can also choose to directly replace $A_6$ by $A_1$ and $A_4$ although it does not reduce the number of vertices, since we now know from theorem \ref{t:a4red} that any graph which includes $A_4$ is reducible.

This leaves us with only the atoms
$A_1$, $A_2$, $A_3$ and $A_7$ for which we do not have reductions, and $A_2$ and $A_3$ must have their vertices of degree 3 adjacent to pendant triangles, which are copies of $A_1$, of course. Similarly, no vertices of degree 2 in an $A_1$ or $A_7$ can be adjacent to common neighbours, as we saw in section \ref{s:irred}.
Moreover, the triangles in $A_7$ are not aloof, and both $A_2$ and $A_3$ need to be adjacent to aloof triangles to be irreducibles, so $A_7$ cannot be in any irreducible graph. 
%Moreover, $A_7$ was shown to be reducible if there was more than one multiple edge in the graph, and so there cannot be more than one copy of $A_7$ as we can use the triple edge in one copy to reduce the other one using lemma \ref{l:2s}, forming an $A_6$. 

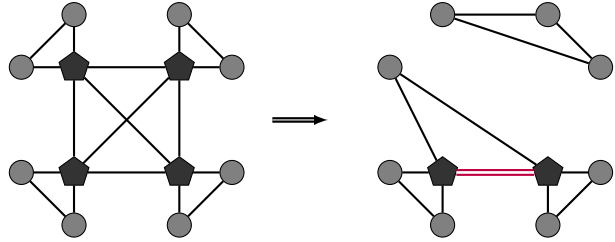
\begin{figure}[ht]
\begin{center}
\begin{tikzpicture}[scale=0.7, transform shape]
        \begin{scope}[shift={(0,0)}]
\foreach \pos/\name in {{(0,0)/a},{(2,0)/b},{(0,2)/c},{(2,2)/d}}
        \node[pnode] (\name) at \pos {};
 \foreach \pos/\name in {{(-1,0)/a1},{(3,0)/b1},{(-1,2)/c1},{(3,2)/d1},{(0,-1)/a2},{(2,-1)/b2},{(0,3)/c2},{(2,3)/d2}}
        \node[enode] (\name) at \pos {};
   \foreach \source/\dest in {a/b,b/c,a/c,a/d,b/d,c/d,a/a1,a/a2,a1/a2,b/b1,b/b2,b1/b2,c/c1,c/c2,c1/c2,d/d1,d/d2,d1/d2}
        \path[edge] (\source) -- (\dest);
 \end{scope}
         \begin{scope}[shift={(7,0)}]
         \node[math,left] () at (-2,1) {{$\tikzLongrightarrow$}};
\foreach \pos/\name in {{(0,0)/a},{(2,0)/b}}
        \node[pnode] (\name) at \pos {};
 \foreach \pos/\name in {{(-1,0)/a1},{(3,0)/b1},{(-1,2)/c1},{(3,2)/d1},{(0,-1)/a2},{(2,-1)/b2},{(0,3)/c2},{(2,3)/d2}}
        \node[enode] (\name) at \pos {};
   \foreach \source/\dest in {a/a1,a/a2,a1/a2,b/b1,b/b2,b1/b2,b/c1,c1/a,c2/d1,c2/d2,d1/d2}
        \path[edge] (\source) -- (\dest);
          \foreach \source/\dest in {a/b}
        \path[dedge] (\source) -- (\dest);
        \end{scope}
\end{tikzpicture}
\end{center}
\caption{Reduction of $A_2$ using $A_1$ and $A_3$}\label{f:a2red}
\end{figure}

Because of the necessary pendant triangles around $A_3$, it is now possible to reduce $A_2$ using it as shown in figure \ref{f:a2red}. We delete all edges in $A_2$, as well as the edges of the aloof triangles adjacent to it, we have 8 vertices which need 2 more edges added to them in order to get a quintic graph. We can accomplish this using $A_3$, its aloof triangles and $A_1$, which together have 8 vertices of degree 2, and only 2 vertices of degree 5, compared to the previous 4. 
 
%However, this fails for the graph in figure \ref{f:k34} since the triangle added needs to be between ... it actually works, but makes double edges.

Recall that the graphs in figures \ref{f:fourv} and \ref{f:ircut} could not be reduced using our given operations. These four small graphs will be our base set of 5-regular graphs with the triangle property. There are no other such graphs with 4 vertices and all others with 6 vertices have been shown to contain at least one of the atoms 
%from table \ref{t:atoms23} 
and we have proven that all of these atoms can be reduced apart from $A_1$ and $A_3$.

\begin{thm}\label{t:class}
All foundational connected quintic graphs with the triangle property and at least eight vertices are constructed from a line graph of a cubic graph ($H$), with a perfect matching $M$, by adding a second edge to $H$ for every edge in $M$.
%either contain one of the given subgraphs or be formed from a $(3,4)$-biregular graph.
\end{thm}

\begin{proof}
Let $G$ be a quintic graph with the triangle property and at least eight vertices, so it is not part of the base set.
From the previous reductions, the only subgraphs of $G$ which can't be reduced by the previous operations are $A_1$ and $A_3$, so there are no diamonds or $K_4$s in $G$, only triangles. At each vertex there must be exactly one double edge and so on removal of the double edges we get $H$ and $M$ as described, and, from \cite{ar:PR}, this is a line graph of a cubic graph.

Note that if we have $t$ vertex disjoint copies of $A_1$ and $s$ copies of $A_3$ they need to be combined in an edge-disjoint way to give vertices of degree 5, so we have $2s$ vertices that start with degree 3, and $3t+s$ with degree 2, and there need to be the same number of each, so $2s=3t+s$, or $s=3t = |V(G)|$.
\end{proof}

\begin{figure}[ht]
\begin{center}
\begin{tikzpicture}[scale=0.5, transform shape]
        \begin{scope}[shift={(0,0)},rotate=30]
\def\ra{.7}
\def\rb{3}
      \draw \foreach \x in {0,120,240}
    {
	   (\x:\ra) -- (\x+120:\ra)
         (\x+60:\rb) -- (\x+180:\rb)
         (\x+60:\rb) -- (\x:\ra)
         (\x-60:\rb) -- (\x:\ra)
          };
      \path[dedge] (60:\rb) -- (180:\rb); 
  \path[dedge] (120:\ra) -- (240:\ra);
    \path[dedge] (300:\rb) -- (0:\ra);
     \foreach \x in {0,120,240} {\node[enode]  at (\x:\ra) {};\node[enode]  at (\x+60:\rb) {};}
 \end{scope}
\node[math,left,color=brown] () at (5.8,0.9) {{\LARGE T}};
\node[math] () at (11.2,.9) {\Large $\tikzLongrightarrow$};
 \begin{scope}[shift={(7.5,0.9)},rotate=45] 
 \def\ra{.7}
\def\rb{1.5}
\def\rc{3.3}
      \draw \foreach \x in {0,90,180,270}
    {
	   (\x:\ra) -- (\x+90:\ra)
         (\x+45:\rb) -- (\x:\ra)
         (\x-45:\rb) -- (\x:\ra)
         (\x+45:\rb) -- (\x:\rc)
         (\x-45:\rb) -- (\x:\rc)
         (\x+90:\rc) -- (\x:\rc)
           };
   \foreach \x in {0,90,180,270} {\path[dedge] (\x+45:\rb) -- (\x:\ra);}
  \path[dedge] (0:\rc) -- (90:\rc); 
  \path[dedge] (180:\rc) -- (270:\rc);
    \foreach \x in {0,90,180,270} {\node[enode]  at (\x:\ra) {};\node[enode]  at (\x:\rc) {};\node[enode]  at (\x+45:\rb) {};}
 \end{scope}
     \begin{scope}[shift={(15,0.9)},rotate=45]
\def\ra{.7}
\def\rb{1.5}
\def\rc{3.3}
      \draw \foreach \x in {0,90,180,270}
    {
	   (\x+90:\rc) -- (\x:\rc)
           };
 
%           \foreach \source/\dest in {a/b}
%        \path[dedge] (\source) -- (\dest);
    \path[dedge] (180:\rc) -- (90:\rc); 
    \path[dedge] (45:\rb) -- (0:\ra); 
    \path[dedge] (-45:\rb) -- (-90:\ra); 
     \foreach \x in {0,180,270} {\path[ledge] (\x+45:\rb) -- (\x:\rc);\path[ledge] (\x+45:\rb) -- (\x+90:\rc);}
\path[ledge] (315:\rb) -- (0:\ra);
\path[ledge] (270:\ra) -- (0:\ra);
\path[ledge] (270:\ra) -- (225:\rb);
     \path[ledge] (135:\ra) -- (0:\ra);
     \path[ledge] (135:\ra) -- (270:\ra);
    \path[ledge] (135:\ra) -- (45:\rb); 
    \path[dedge] (225:\rb) -- (135:\ra);
     \path[ledge] (0:\rc) .. controls (-85:5) and (-95:5) .. (180:\rc);
     \path[ledge] (90:\rc) .. controls (5:5) and (-5:5) .. (270:\rc);
  \node[enode,color=brown]  at (135:\ra) {};   
    \foreach \x in {0,90,180,270} {\node[enode]  at (\x:\rc) {};}
      \foreach \x in {0,180,270} {\node[enode]  at (\x+45:\rb) {};   }
           \foreach \x in {0,270} {\node[enode]  at (\x:\ra) {}; }
 \end{scope}
\end{tikzpicture}
\end{center}
\caption{Small foundational graphs}\label{f:final}
\end{figure}
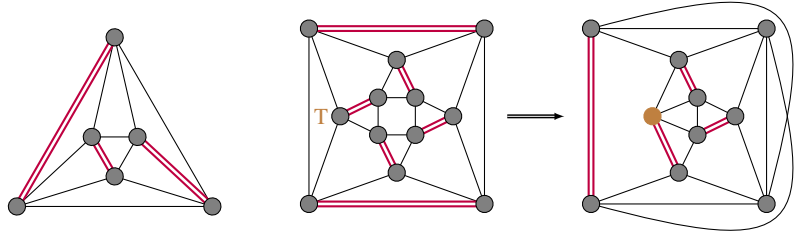

The first couple of graphs in this family are shown in figure \ref{f:final}, from $K_4$ and the cube. If a cubic graph has $2n$ vertices then it has $\frac{3n}{2}$ edges and that is the number of vertices in the quintic graph formed from the line graph; that number has to be even to form a perfect matching, so only cubic graphs with a number of vertices that is a multiple of 4 can be used. 

Note that there is a way to reduce the second graph, as shown, by focusing on a triangle which has no multiple edges. such as $T$. We can delete one vertex $v$ of $T$ and double the other edge of $T$. It is possible to contract $v$'s neighbours outside of $T$ as they are part of an aloof triangle, and then we join the other vertices of $T$ to the vertex of each other's pendant triangles that was originally part of a double edge.
The graph that remains is quintic and has the triangle property still. This cannot be done for the graph with six vertices as the pendant triangles are not vertex disjoint.

\subsection{Future implications}

We have shown how to reduce all quintic graphs with the triangle property, to the specific families. It is also possible to reverse the reductions and generate larger graphs. In particular, both the Z-reduction and X-reductions have a very simple application; for the former we can pick two non-adjacent vertices and try to split the neighbours of each into sets of 2 and 3 so that the triangle property is preserved.

It is possible to proceed similarly for two adjacent vertices with the X-reductions, and, experimentally, most graphs with the triangle property can be generated in this way, but not all. Additionally, both basic operations only introduce simple edges, so if we start with a simple graph, a larger simple graph will be formed. Unfortunately, not all simple quintic graphs with the triangle property can be constructed this way.

Note also that the Z-reduction preserves planarity, so can be used to create arbitrarily large quintic planar graphs, as well as those with the triangle property. However, there are some simple quintic planar graphs which do not arise in this way such as the unique planar quintic graph with 16 vertices and diameter 3 from \cite{ar:diam3}; more operations are required as in \cite{ar:hmr}.

The middle graph in figure \ref{f:final} is from the line graph of the cube and, since the cube is triangle-free, there are exactly 8 triangles in the 5-regular graph with 12 vertices formed from it. This isn't a contradiction to the result in \cite{ar:jp} since in that article triangles are counted by edges rather than vertices, with multiple edges giving multiple triangles, so there are 16 triangles considering edges in the right hand graph under that criterion.
% Not all quartic graphs have perfect matchings, 
% https://cfwebprod.sandia.gov/cfdocs/CompResearch/docs/graph.pdf

\end{document}